\documentclass[a4paper,10pt]{article}

\usepackage{amsthm}

\usepackage{amsmath}

\usepackage{amsfonts}

\usepackage{amssymb}

\usepackage{amsmath,amsthm,amscd,amssymb}

\usepackage{latexsym}

\usepackage{graphicx}

\usepackage[all]{xy}

\theoremstyle{plain}

\newtheorem{thm}{\textbf{Theorem}}

\newtheorem{cn}{\textbf{Conjecture}}

\newtheorem{lem}{\textbf{Lemma}}

\newtheorem{df}{\textbf{Definition}}

\newtheorem{cor}{\textbf{Corollary}}

\newtheorem{prop}{\textbf{Proposition}}

\newcommand{\A}{\Bbb{A}}

\newcommand{\R}{\Bbb{R}}

\newcommand{\C}{\Bbb{C}}

\newcommand{\X}{\Bbb{X}}

\newcommand{\Q}{\Bbb{Q}}

\newcommand{\F}{\Bbb{F}}

\newcommand{\G}{\Bbb{G}}

\newcommand{\Z}{\Bbb{Z}}

\newcommand{\T}{\Bbb{T}}

\newcommand{\p}{\frak{p}}

\newcommand{\m}{\frak{m}}

\newcommand{\K}{\mathcal{K}}

\newcommand{\Hom}{\text{Hom}}

\newcommand{\End}{\text{End}}

\newcommand{\tr}{\text{tr}}

\newcommand{\im}{\text{im}}

\newcommand{\Gal}{\text{Gal}}

\newcommand{\GL}{\text{GL}}

\newcommand{\y}{\hspace{6pt}}

\title{{\bf{Eigenvarieties and invariant norms: Towards $p$-adic Langlands for $U(n)$}}}

\author{Claus M. Sorensen}

\begin{document}

\date{Preprint}

\maketitle

\begin{abstract} 
We give a proof of the Breuil-Schneider conjecture in a large number of cases, which {\it{complement}} the indecomposable case, which we dealt with earlier in [Sor]. In some sense, only the Steinberg representation lies at the intersection of the two approaches. In this paper, we view the conjecture from a broader global perspective. If $U_{/F}$ is any definite unitary group, which is an inner form of $\GL(n)$ over $\K$, we point out how the eigenvariety $\X(K^p)$ parametrizes a global $p$-adic Langlands correspondence between certain $n$-dimensional $p$-adic semisimple representations $\rho$ of $\Gal(\bar{\Q}|\K)$ (or what amounts to the same, pseudo-representations) and certain Banach-Hecke modules $\mathcal{B}$ with an admissible unitary action of $U(F\otimes \Q_p)$, when $p$ splits. We express the locally regular-algebraic vectors of $\mathcal{B}$ in terms of the Breuil-Schneider representation of $\rho$. Upon completion, this produces a {\it{candidate}} for the $p$-adic local Langlands correspondence in this context.
As an application, we give a weak form of local-global compatibility in the crystalline case, showing that the Banach space representations $B_{\xi,\zeta}$ of Schneider-Teitelbaum [ScTe] fit the picture as predicted. There is a compatible global mod $p$ (semisimple) Langlands correspondence parametrized by $\X(K^p)$. We introduce a natural notion of refined Serre weights, and link them to the existence of crystalline lifts of prescribed Hodge type and Frobenius eigenvalues. At the end, we give a rough candidate for a local mod $p$ correspondence, formulate a local-global compatibility conjecture, and explain how it implies the conjectural Ihara lemma in [CHT].

\footnote{{\it{2000 AMS Mathematics Classification}}: 11F33.}
\end{abstract}


\section{Introduction}

Let $\K$ be a number field. The Fontaine-Mazur conjecture [FM] predicts a characterization of all (irreducible) Galois representations $\rho:\Gamma_{\K}=\Gal(\bar{\Q}|\K)\rightarrow \GL_n(\Q_p)$ occurring naturally. That is, in the etale cohomology $H^{\bullet}(X,\Q_p)$ of some smooth projective variety $X_{/\K}$. It is a major result (due to Tsuji and others) that every such $\rho$ is {\it{geometric}}, which means it is unramified at all but finitely many places, and potentially semistable at all places above $p$. Fontaine and Mazur assert the converse, that every geometric $\rho$ occurs in cohomology (up to a Tate twist). The potentially semistable representations are now more or less completely understood, by work of Colmez and Fontaine [CM]. They are given by admissibly filtered $(\phi,N)$-modules (with Galois action), which are objects of a more concrete combinatorial nature. The $p$-adic Langlands program, still in its initial stages, attempts to link $p$-adic Hodge theory with non-archimedean functional analysis. Locally, if $K$ is a fixed finite extension of $\Q_p$, and $L|\Q_p$ is another sufficiently large finite extension (the coefficient field), one hopes to pair certain Galois representations $\rho:\Gamma_K \rightarrow \GL_n(L)$ with certain Banach $L$-spaces with a unitary admissible $\GL_n(K)$-action. This is now well-understood for $\GL_2(\Q_p)$ thanks to recent work of Berger, Breuil, Colmez, Paskunas, and others. See [Bg] for a nice survey. The goal of this paper is to shed some light on a {\it{global}} analogue, for any $n$, and any CM field $\K$. To give the flavor, if $\K|\Q$ is a quadratic imaginary field, in which $p$ splits, we will set up a bijection between certain Galois representations $\rho:\Gamma_{\K}\rightarrow \GL_n(L)$ (actually, pseudo-representations) and certain Banach-Hecke modules with a unitary admissible $\GL_n(\Q_p)$-action. This is most likely folklore. More importantly, we relate the algebraic vectors to the $p$-adic Hodge theory on the Galois side. The word {\it{certain}} here has a precise meaning. It means those representations which {\it{come from an eigenvariety}}, of some fixed tame level $K^p$. We will be precise below.

\medskip

\noindent We model the discussion on the $\GL(2)$-case: With any continuous Galois representation $\rho:\Gamma_{\Q_p}\rightarrow \GL_2(L)$, the $p$-adic Langlands correspondence associates a unitary Banach $L$-space representation $B(\rho)$ of $\GL_2(\Q_p)$. Moreover, $\rho$ is de Rham with distinct Hodge-Tate weights precisely when $B(\rho)^{alg}\neq 0$. If so, let us say the weights are $\{0,1-k\}$ with $k \geq 2$ (with the convention that the cyclotomic character has weight $-1$), then the algebraic vectors are given by
$$
B(\rho)^{alg}=\text{Sym}^{k-2}(L^2)\otimes \pi(\rho)
$$
for a smooth generic representation $\pi(\rho)$, possibly reducible, obtained by a slight modification of the classical local Langlands correspondence.

\medskip

\noindent {\it{The Breuil-Schneider conjecture}}. The local $p$-adic Langlands program is somewhat vague, and a precise conjectural framework is still developing, beyond the case of $\GL_2(\Q_p)$, where pretty much everything is known. However, there is a weak (but precise) almost "skeletal" version formulated in [BrSc], which we now recall. We keep our finite extension $K|\Q_p$, and a finite Galois extension thereof $K'|K$. Pick a third field of coefficients $L\subset \bar{\Q}_p$, finite over $\Q_p$, but large enough so that it contains the Galois closures of $K$ and $K_0'$ (the maximal unramified subfield of $K'$). 
The roles of these fields are the following. We consider potentially semistable representations $\rho:\Gamma_K \rightarrow \GL_n(L)$, which becomes semistable when restricted
to $\Gamma_{K'}$. As mentioned above, such $\rho$ correspond to $(\phi,N)\times \Gal(K'|K)$-modules $D$, with an admissible filtration. This makes use of Fontaine's period ring $B_{st}$,
$$
D=(B_{st}\otimes_{\Q_p}\rho)^{\Gamma_{K'}}.
$$ 
This is a finite free $K_0'\otimes_{\Q_p} L$-module, of rank $n$, with a semilinear Frobenius $\phi$, a (nilpotent) monodromy operator $N$ such that $N\phi=p\phi N$, a commuting action of $\Gal(K'|K)$, and an admissible Galois-stable filtration on $D_{K'}$. Note 
$$
K'\otimes_{\Q_p}L\simeq {\prod}_{\tau\in \Hom(K,L)}K'\otimes_{K,\tau}L.
$$
Accordingly, $D_{K'}\simeq \prod_{\tau}D_{K',\tau}$, and each $K'\otimes_{K,\tau}L$-module $D_{K',\tau}$ is filtered.

\begin{itemize}
\item {\it{Hodge-Tate numbers}}. For each $\tau:K \hookrightarrow L$, we let $i_{1,\tau}\leq \cdots\leq i_{n,\tau}$ denote the jumps in the Hodge filtration (listed with multiplicity). That is,
$$
\text{gr}^i(D_{K',\tau})\neq 0 \Leftrightarrow i\in \{i_{1,\tau},\ldots, i_{n,\tau}\}.
$$
We will denote this multiset of integers by $HT_{\tau}(\rho)=\{i_{j,\tau}:j=1,\ldots,n\}$.
\item {\it{Weil-Deligne representation}}. If we forget about the filtration, the resulting $(\phi,N)\times \Gal(K'|K)$-module corresponds to a Weil-Deligne representation, once we fix an embedding $K_0'\hookrightarrow L$. This goes as follows, see Proposition 4.1 in [BrSc]. The underlying $n$-dimensional $L$-vector space is
$$
D_L=D\otimes_{K_0'\otimes_{\Q_p} L}L,
$$
with the induced $N$ coming from $B_{st}$, and with $r:W_K \rightarrow \GL(D_L)$ defined by
$r(w)=\phi^{-d(w)}\circ \bar{w}$. Here $\bar{w}$ denotes the image of $w$ in $\Gal(K'|K)$, and $d(w)$ gives the power of arithmetic Frobenius which $w$ induces. The ensuing Weil-Deligne representation becomes unramified upon restriction to $W_{K'}$. We will denote it by $WD(\rho)=(r,N,D_L)$ throughout the text.
\end{itemize}

\noindent The Breuil-Schneider conjecture asks for a characterization of the data arising in this fashion, assuming all Hodge-Tate numbers are distinct. To state it, start with abstract data. Firstly, for each embedding $\tau:K \hookrightarrow L$, say we are given $n$ distinct integers $HT_{\tau}=\{i_{1,\tau}<\cdots<i_{n,\tau}\}$. Secondly, say we are given some $n$-dimensional Weil-Deligne representation $WD$, with coefficients in $L$, which become unramified on $W_{K'}$. With these data, below we will associate a locally algebraic representation $BS$ of $\GL_n(K)$, with coefficients in $L$. The algebraic part is defined in terms of the $HT_{\tau}$, the smooth part in terms of $WD$. Our data should come from a $\rho$ precisely when $BS$ has a $\GL_n(K)$-stable $\mathcal{O}_L$-lattice; the unit ball of an invariant norm. The following is Conjecture 4.3 in [BrSc], also announced as Conjecture 4.1 in Breuil's 2010 ICM address [Bre]:

\medskip

\noindent {\bf{The Breuil-Schneider conjecture}}. {\it{The following are equivalent:
\begin{itemize}
\item[(1)] The data $HT_{\tau}$ and $WD^{F-ss}$ arise from a potentially semistable $\rho$.
\item[(2)] $BS$ admits a norm $\|\cdot\|$, invariant under the action of $\GL_n(K)$.
\end{itemize}
}}

\medskip

\noindent Before we recall the status of the conjecture, we return to the definition of $BS$.

\begin{itemize}
\item {\it{Algebraic part}}. Introduce $b_{j,\tau}=-i_{n+1-j,\tau}-(j-1)$. That is, write the $i_{j,\tau}$ in the opposite order, change signs, and subtract $(0,1,\ldots,n-1)$.
We let $\xi_{\tau}$ be the irreducible algebraic $L$-representation of $\GL_n$, of highest weight 
$$
b_{1,\tau}\leq b_{2,\tau}\leq \cdots \leq b_{n,\tau}
$$
relative to the {\it{lower}} triangular Borel. Their tensor product $\xi=\otimes_{\tau}\xi_{\tau}$, with $\tau$ running over $\Hom(K,L)$, is then an irreducible algebraic representation of $\GL_n(K\otimes_{\Q_p}L)$ over $L$. We will view $\xi$ as a representation of $\GL_n(K)$.

\item {\it{Smooth part}}. By the classical local Langlands correspondence [HT], the Frobenius-semisimplification $WD^{F-ss}\simeq \text{rec}_n(\pi^{\circ})$, for some irreducible admissible smooth representation $\pi^{\circ}$ of $\GL_n(K)$, defined over $\bar{\Q}_p$. Here $\text{rec}_n$ is normalized as in [HT]. To define it over $\bar{\Q}_p$, we need to fix a square-root $q^{\frac{1}{2}}$, where $q=\#\F_K$. By the Langlands classification, one has
$$
\text{Ind}_P(Q(\Delta_1)\otimes \cdots \otimes Q(\Delta_r))\overset{!}{\twoheadrightarrow} \pi^{\circ},
$$
a unique irreducible quotient, where the $Q(\Delta_i)$ are generalized Steinberg representation built from the $\Delta_i$, which are segments of supercuspidals, suitably ordered.
The smooth part of $BS$ is now defined to be
$$
\pi=\text{Ind}_P(Q(\Delta_1)\otimes \cdots \otimes Q(\Delta_r))\otimes |\det|^{\frac{1-n}{2}},
$$
or rather its model over $L$, which is independent of the choice of $q^{\frac{1}{2}}$. Note
that $\pi\simeq \pi^{\circ}\otimes |\det|^{(1-n)/2}$ if and only if $\pi^{\circ}$ is generic (that is, has a Whittaker model). For that reason, the association $WD \mapsto \pi$ is often called the {\it{generic}} local Langlands correspondence.
\end{itemize}

\noindent We let $BS=\xi\otimes_L \pi$, following [BrSc] (although they do not use the notation $BS$). In fact, we will find it more convenient to work with a different normalization. In the above construction there is a choice of a {\it{sign}}; essentially reflected in whether one twists by $|\det|^{(1-n)/2}$ or its inverse. The latter is more commonly used in the references we rely on. The resulting representation is just a twist of $BS$ by a harmless explicit $p$-adically unitary continuous character. Namely,
$$
\text{$\widetilde{BS}=BS \otimes_L \mu^{n-1}$, $\y$ $\mu(g)=N_{K|\Q_p}(\det g)^{\times}$,}
$$
where $a^{\times}=a|a|_p=BS(\chi_{cyc})(a)\in \Z_p^{\times}$ denotes the unit factor of an $a \in \Q_p^*$. Of course, $\widetilde{BS}$ has an invariant norm if and only if $BS$ does, so it makes no real difference. It reflects a Tate twist: $\widetilde{BS}(\rho)$ is nothing but $BS(\rho\otimes \chi_{cyc}^{n-1})$.

\medskip

\noindent The implication (2) $\Rightarrow$ (1) in the conjecture is in fact completely known. After many cases were worked out in [ScTe] and [BrSc], the general case was settled by Y. Hu in his thesis [Hu]. In fact Hu proves that (1) is equivalent to the {\it{Emerton condition}}, which is a purely group-theoretical condition:
$$
\text{(3)      $J_P(BS)^{Z_M^+=\chi}\neq 0 \Longrightarrow \forall z \in Z_M^+: |\delta_P^{-1}(z)\chi(z)|_p\leq 1$.}
$$
Here $J_P$ is Emerton's generalization of the Jacquet functor, introduced and studied in [Em1] and [Em2]. The heart of Hu's proof is to translate (3) into finitely many inequalities relating the Hodge polygon to the Newton polygon. In the vein [FoRa], he is then able to show the existence of an admissible filtration compatible with the given data.
The implication (2) $\Rightarrow$ (3) is relatively easy.

\medskip

\noindent What remains, is to produce an invariant norm on $BS(\rho)$, for any potentially semistable $\rho$ (with distinct Hodge-Tate weights). 
One of the main motivations for writing this paper, was to make progress in this direction, (1) $\Rightarrow$ (2). We proved this in [Sor] when $WD(\rho)$ is indecomposable
(in other words, $\pi^{\circ}=Q(\Delta)$ is generalized Steinberg). Here the Emerton condition boils down to just integrality of the central character, and in fact the resulting conjecture was stated explicitly as 5.5 in [BrSc]. The key point of [Sor] was to make use of the fact that $Q(\Delta)$ is a discrete series representation, and therefore admits a {\it{pseudo-coeffcient}}. Inserting this as a test-function in the trace formula for a certain definite unitary group, one can pass to a global setup (a la Grunwald-Wang). Finally, the desired norm was found by relating classical $p$-adic algebraic modular forms to completed cohomology $\tilde{H}^0$, as introduced in great generality in [Emer]. This whole argument is purely group-theoretical, and in fact carries over to any connected reductive group over $\Q_p$, exploiting a compact form (using a Galois cohomological computation of Borel and Harder, which shows the existence of locally prescribed forms). We should point out that the supercuspidal case is much easier. In this case there are several ways to produce a norm (compact induction, for example). 

\medskip

\noindent One of the outcomes of this paper, is a {\it{complement}} to the main result of [Sor]. The idea of relating algebraic modular forms to $\tilde{H}^0$, already present in [Emer], can be pushed further, now that local-global compatibility at $p=\ell$ is available in the "book project" context. This was proved recently by Barnet-Lamb, Gee, Geraghty, and Taylor in the so-called Shin-regular case [BGGT], and this regularity hypothesis was then shown to be unnecessary by Caraiani, as part of her 2012 Harvard Ph.D. thesis [Car]. This results in the following somewhat vague Theorem A, which we will make more precise in Theorem B below.
\medskip

\noindent {\bf{Theorem A}}. {\it{The Breuil-Schneider conjecture holds for potentially semistable $\rho$, which come from a regular, classical, irreducible point on a unitary eigenvariety. }}

\medskip

\noindent {\it{Eigenvarieties}}. We will combine the approaches of [Chen] and [Emer]. Thus let $\K$ be a CM field, with maximal totally real subfield $F$. Let $D$ be a central simple $\K$-algebra of $\dim_{\K}(D)=n^2$, equipped with an anti-involution $\star$ of the second kind (that is, $\star|_{\K}$ is conjugation). We introduce the unitary $F$-group $U=U(D,\star)$, an outer form of $\GL(n)$, which becomes the inner form $D^{\times}$ over $\K$. It will be convenient to also introduce $G=\text{Res}_{F|\Q}(U)$. We will always assume $U$ is totally definite. In other words, that $G(\R)$ is a {\it{compact}} Lie group,
therefore a product of copies of $U(n)$ (hence also connected).

\medskip

\noindent We will fix a prime number $p$ such that every place $v|p$ of $F$ splits in $\K$, and such that $D_w^{\times}\simeq \GL_n(\K_w)$ for every $w|v$. To keep track of various identifications, it is customary to {\it{choose}} a place $\tilde{v}$ of $\K$ above every $v|p$. Once and for all, also fix
an isomorphism $\iota: \C \overset{\sim}{\longrightarrow} \bar{\Q}_p$. This gives rise to an identification
$$
\Hom(F,\R)=\Hom(F,\C)\simeq \Hom(F,\bar{\Q}_p)=\sqcup_{v|p}\Hom(F_v,\bar{\Q}_p),
$$
and similarly for $\Hom(\K,\C)$. By assumption $F_v \simeq \K_w$ for $w|v$, so the choices $\{\tilde{v}\}$ just amount to fixing a CM-type $\Phi$, which is ordinary for $\iota$, in the sense of [Katz]. This will ensure that the various identifications we make are compatible.

\medskip

\noindent The eigenvariety for $G$ depends on the choice of a tame level $K^p \subset G(\A_f^p)$. It is a reduced rigid analytic space $\X_{/E}$, where we take $E$ to be the Galois closure of $F$ in $\bar{\Q}_p$, with additional structure: 
$$
\text{$\chi:\X \rightarrow \hat{T}$, $\y$ $\lambda: \mathcal{H}(K^p)^{\text{sph}}\rightarrow \mathcal{O}(\X)$.}
$$
Here $\hat{T}_{/E}$ is weight space, parametrizing locally analytic characters of $T(\Q_p)$, and $\mathcal{H}(K^p)^{\text{sph}}$ is the spherical central subalgebra of the Hecke $E$-algebra $\mathcal{H}(K^p)$. Moreover, there is a Zariski-dense subset $X_{cl}\subset \X(\bar{\Q}_p)$ such that the evaluation
$$
\text{$\X(\bar{\Q}_p)\longrightarrow (\hat{T}\times \text{Spec}\mathcal{H}(K^p)^{\text{sph}})(\bar{\Q}_p)$, $\y$ $x \mapsto (\chi_x,\lambda_x)$,}
$$
identifies $X_{cl}$ with the set of {\it{classical}} points. Roughly this means that, first of all $\chi_x=\psi_x\theta_x$ is locally algebraic, and there exists an automorphic representation $\pi$ of weight $\psi_x$ such that $\pi_p\hookrightarrow \text{Ind}_{B}^{G}(\theta_x)$, and $\mathcal{H}(K^p)^{\text{sph}}$ acts on $\pi_f^{K^p}\neq 0$ by the character $\lambda_x$. (The condition that $\pi_p$ embeds in a principal series is the analogue of the "finite slope" requirement showing up in the classical works of Coleman, Mazur and others. The choice of a $\theta_x$ is called a refinement of $\pi$.)

\medskip

\noindent It is of utmost importance to us that the eigenvariety carries a family of Galois representations. To be more precise, if we let $\Sigma=\Sigma(K^p)$ be the set of ramified places, there is a unique continuous $n$-dimensional pseudo-representation 
$$
\mathcal{T}: \Gamma_{\K,\Sigma}\rightarrow \mathcal{O}(\X)^{\leq 1}
$$
associated with $\lambda: \mathcal{H}(K^p)^{\text{sph}}\rightarrow \mathcal{O}(\X)$, in the sense that for all places $w \notin \Sigma$,
$$
\mathcal{T}(\text{Frob}_w)=\lambda(b_{w|v}(h_w)).
$$
Here $h_w$ is the element of the spherical Hecke algebra for $\GL_n(\K_w)$, which acts via the sum of the (integral) Satake parameters on spherical vectors, and $b_{w|v}$ is the standard base change homomorphism between the pertaining spherical Hecke algebras, see [Min]. In particular, by Procesi and Taylor, for each $x \in \X(\bar{\Q}_p)$ there is a unique semisimple Galois representation 
$$
\text{$\rho_x:\Gamma_{\K,\Sigma}\rightarrow \GL_n(\bar{\Q}_p)$, $\y$ $\mathcal{T}_x=\tr(\rho_x)$.}
$$
In fact, the way $\mathcal{T}$ is constructed, is by first defining $\rho_x$ for {\it{regular}} classical points $x\in X_{reg}$, by which we mean the dominant character $\psi_x$ is given by a strictly decreasing sequence of integers (at some place). Thanks to [Whi], this guarantees that $\pi$ has a base change to $\GL_n(\A_{\K})$ of the form $\Pi=\boxplus \Pi_i$, where the $\Pi_i$ are cohomological {\it{cuspidal}} (as opposed to just discrete) automorphic representations to which one can attach Galois representations. Now, $X_{reg}$ can be shown to be Zariski-dense, and a formal argument in [Che] interpolates the pseudo-characters $\tr(\rho_x)$ for $x\in X_{reg}$ by a unique $\mathcal{T}$, which one can then specialize at {\it{any}}
point $x \in \X(\bar{\Q}_p)$.

\medskip

\noindent Now we can clarify the statement in Theorem A: If $x \in X_{reg}$ is a classical point such that $\rho_x$ is irreducible (globally, as a representation of $\Gamma_{\K}$), and $w|p$ is a place of $\K$, then $\rho_x|_{\Gamma_{\K_w}}$ is potentially semistable, {\it{and}} its locally algebraic representation $BS(\rho_x|_{\Gamma_{\K_w}})$ admits a $\GL_n(\K_w)$-invariant norm $\|\cdot\|$.

\medskip

\noindent {\it{A global $p$-adic Langlands correspondence}}. The actual construction of an invariant norm $\|\cdot\|$ is more interesting than its mere existence. It comes out of a much more precise result, which we now describe. Fix a finite extension $L|E$. At each point $x \in \X(L)$, we have pseudo-representation $\mathcal{T}_x:\Gamma_{\K,\Sigma}
\rightarrow L$ (the trace of an actual representation $\rho_x$, which may or may not be defined over $L$). On the other hand, with $x \in \X(L)$ we associate the Banach $L$-space
$$
\mathcal{B}_x=(L\otimes_E \tilde{H}^0(K^p))^{\frak{h}=\lambda_x}
$$
where $\frak{h}=\mathcal{H}(K^p)^{\text{sph}}$ is shorthand notation. This space is really very concrete. The completed cohomology $\tilde{H}^0(K^p)$, defined in much greater generality by Emerton, is here nothing but the space of all {\it{continuous}} functions
$$
\text{$f: Y(K^p)\rightarrow E$, $\y$ $Y(K^p)=\underset{K_p}{\varprojlim} Y(K_pK^P)$, $\y$ $Y(K)=G(\Q)\backslash G(\A_f)/K$,}
$$
with supremum norm. The superscript $\frak{h}=\lambda_x$ means we take the eigenspace for the character $\lambda_x: \frak{h}\rightarrow L$ (not the generalized eigenspace). Note that $\mathcal{B}_x$ is much more than just a Banach $L$-space: For one thing, it is a Banach module for the Banach-Hecke algebra $\hat{\mathcal{H}}(K^p)$ (see [ScTe] for a detailed discussion of these). For another thing, there is a natural $\hat{\mathcal{H}}(K^p)$-linear action of $G(\Q_p)$ by {\it{isometries}} of $\mathcal{B}_x$, which is admissible (meaning that its reduction $\bar{\mathcal{B}}_x$ is a smooth admissible representation of $G(\Q_p)$ over $\F_L$, in the usual sense). Now, for $x,x'\in \X(L)$,
$$
\mathcal{T}_x=\mathcal{T}_{x'} \Leftrightarrow \lambda_x=\lambda_{x'}\Leftrightarrow \mathcal{B}_x=\mathcal{B}_{x'},
$$
since each $b_{w|v}$ is onto, see Corollary 4.2 in [Min] (a fact also used on p. 10 of [CHL]). In other words, the set of all pairs $(\mathcal{T}_x,\mathcal{B}_x)$, with $x \in \X(L)$, is the graph of a bijection between the images of the two projections:
$$
\left\{ \begin{matrix} \text{$n$-dimensional pseudo-representations} \\ \text{$\mathcal{T}: \Gamma_{\K,\Sigma}\rightarrow L$ coming from $\X(L)$} \end{matrix} \right\}
\longleftrightarrow
$$
$$
\left\{ \begin{matrix} \text{Banach $\hat{\mathcal{H}}_L(K^p)$-modules $\mathcal{B}$ with admissible} \\ \text{unitary $G(\Q_p)$-action, coming from $\X(L)$}  \end{matrix} \right\}.
$$
Here $\mathcal{T}\leftrightarrow \mathcal{B}$ means there is a point $x \in \X(L)$ such that $\mathcal{T}=\mathcal{T}_x$ and $\mathcal{B}=\mathcal{B}_x$.
(We say that a pseudo-character $\mathcal{T}: \Gamma_{\K,\Sigma}\rightarrow L$ comes from $\X(L)$ if it is of the form $\mathcal{T}_x$ for a point $x\in \X(L)$. Similarly for Banach modules.)

\medskip

\noindent For ease of exposition, let us assume we have {\it{split ramification}}. In other words, $S(K^p)\subset \text{Spl}_{\K|F}$. Then local base change is defined everywhere, and there is a unique automorphic representation $\pi_x$ associated with a point $x\in X_{cl}$ such that $\rho_x$ is irreducible (indeed its global base change is cuspidal and determined almost everywhere). It is expected (and perhaps known?) that $m(\pi_x)=1$. Our main result in this paper is the following, proved in Section 6.

\medskip

\noindent {\bf{Theorem B}}. {\it{Assume $S(K^p)\subset \text{Spl}_{\K|F}$. For each point $x \in X_{reg}\cap X_{irr}$, defined over $L$, such that $m(\pi_x)=1$, there is a unique (up to topological equivalence) 
Banach space $B(\rho_x)$ over $L$ with an \underline{admissible} unitary $G(\Q_p)$-action such that
\begin{itemize}
\item[(1)] $B(\rho_x)^{ralg}\simeq \widetilde{BS}(\rho_x):=\bigotimes_{v|p}\widetilde{BS}(\rho_x|_{\Gamma_{\K_{\tilde{v}}}})$ is \underline{dense}.
\item[(2)] There is a $G(\Q_p)\times \hat{\mathcal{H}}(K^p)$-equivariant topological isomorphism,
$$
B(\rho_x)\otimes (\otimes_{v\nmid p}\pi_{x,v}^{K_v}) \overset{\sim}{\longrightarrow} \overline{\mathcal{B}_x^{ralg}}.
$$ 
(Here $\overline{\mathcal{B}_x^{ralg}}$ denotes the closure of the regular-algebraic vectors in $\mathcal{B}_x$.)
\item[(3)] If $\rho_x$ is crystalline above $p$, there is a continuous $G(\Q_p)$-map, with dense image,
$$
B_{\xi_x,\zeta_x} \longrightarrow B(\rho_x),
$$
which restricts to an isomorphism $H_{\xi_x,\zeta_x} \overset{\sim}{\longrightarrow} B(\rho_x)^{ralg}$. (Here $H_{\xi,\zeta}$ and $B_{\xi,\zeta}$ are the spaces introduced by Schneider and Teitelbaum [ScTe], and we take $\xi_x$ of highest weight $\psi_x$, and $\zeta_x$ to be the eigensystem of $\theta_x$.) 
\end{itemize}
}}

\medskip

\noindent {\it{A local $p$-adic Langlands correspondence?}} The $\GL_2(\Q_p)$-case
suggests there should be a local "correspondence" for any finite extension $K|\Q_p$. Being cautious, there should at least be {\it{map}} $\rho \mapsto \frak{B}(\rho)$ (defined for possibly non-semisimple $\rho$. In other words, for representations, as opposed to just pseudo-representations)
$$
\left\{ \begin{matrix} \text{continuous representations $\rho: \Gamma_K \rightarrow \GL_n(L)$}  \end{matrix} \right\}
\overset{?}{\longrightarrow}
$$
$$
\left\{ \begin{matrix} \text{Banach $L$-spaces $\mathfrak{B}$ endowed with} \\ \text{admissible unitary $\GL_n(K)$-action}  \end{matrix} \right\},
$$
which should map irreducible $\rho$ to topologically irreducible $\frak{B}(\rho)$, and one should be able to recover $\rho$ from $\frak{B}(\rho)$. If so, we believe in the provisional equality
$$
B(\rho_x)\overset{?}{=} \hat{\otimes}_{v|p}\mathfrak{B}(\rho_x|_{\Gamma_{\K_{\tilde{v}}}}),
$$
for $x$ as in Theorem B, at least when the restrictions $\rho_x|_{\Gamma_{\K_{\tilde{v}}}}$ are irreducible. Of course, the flaring question here is whether $B(\rho_x)$ only depends on these restrictions at $p$ (and whether it factors as a tensor product). This is the crux of the matter, and seems hard. We have nothing to offer in this direction here.

\medskip

\noindent For definite unitary groups $U(2)$ over $\Q$, the above question is well-posed: The local $p$-adic Langlands correspondence $\frak{B}$ has been constructed for $\GL_2(\Q_p)$ by Berger, Breuil, Colmez, Kisin, and others, by ingenious use of $(\phi,\Gamma)$-modules; a truly monumental result! In the crystalline case, $\frak{B}$ is (typically) given by the universal modules $B_{\xi,\zeta}$. Thus, at least in this case, Theorem B gives some sort of weak local-global compatibility, in the vein of [Eme]. In joint work with Przemyslaw Chojecki, we are hoping to extends this to the {\it{non}}-crystalline case, by employing the deformation theoretical techniques used in [Eme].

\medskip

\noindent {\it{A global mod $p$ Langlands correspondence}}. In the last section, we introduce the mod $p$ correspondence, for any finite extension $L|E$, and any tame level $K^p$, 
$$
\left\{ \begin{matrix} \text{$n$-dimensional pseudo-representations} \\ \text{$t: \Gamma_{\K,\Sigma}\rightarrow \F_L$ coming from $\X(L)$} \end{matrix} \right\}
\longleftrightarrow
$$
$$
\left\{ \begin{matrix} \text{$\mathcal{H}_{\F_L}(K^p)$-modules $b$ with admissible} \\ \text{$G(\Q_p)$-action, coming from $\X(L)$}  \end{matrix} \right\}.
$$
Here, 
$$
t_x=\tr \bar{\rho}_x^{ss} \longleftrightarrow b_x=H^0(K^p,\F_L)^{\frak{h}^{\circ}=\bar{\lambda}_x}
$$
for any $x \in \X(L)$. If $\omega$ is an irreducible $G(\F_p)$-representation over $\bar{\F}_p$, a Serre weight, $\Hom_{G(\Z_p)}(\omega,b_x)$ is naturally identified with the  
generalized eigenspace for $\bar{\lambda}_x$ in the mod $p$ modular forms $\mathcal{A}_{\omega}(G(\Z_p)K^p,\bar{\F}_p)$. This latter space carries a natural action of $\mathcal{H}_{\omega}(G,K)$, the Hecke algebra at $p$, which is commutative [Her]. This triggers the following notion. We define a {\it{refined}} Serre weight as a pair $(\omega,\nu)$, where $\nu: \mathcal{H}_{\omega}(G,K) \rightarrow \bar{\F}_p$ is an algebra character. Let
$$
(\omega,\nu) \in \mathcal{W}_+(\bar{\rho}_x) \Longleftrightarrow \mathcal{A}_{\omega}(G(\Z_p)K^p,\bar{\F}_p)^{\mathcal{H}_{\omega}(G,K)\otimes\frak{h}^{\circ}=\nu\otimes\bar{\lambda}_x} \neq 0.
$$
Equivalently, $\mathcal{W}_+(\bar{\rho}_x)$ contains $(\omega,\nu)$ when there is a nonzero $G(\Q_p)$-map,
$$
\pi(\omega,\nu):=\text{c-Ind}_{K}^G(\omega)\otimes_{\mathcal{H}_{\omega}(G,K),\nu}\bar{\F}_p \rightarrow b_x.
$$
In Proposition 4 below, we aim at a more Galois-theoretical description of $\mathcal{W}_+(\bar{\rho}_x)$, inspired by [Gee]. Namely, given $(\omega,\nu)$ in there, we note that $\bar{\rho}_x$ has a crystalline lift with Hodge-Tate weights and Frobenius eigenvalues prescribed by  $(\omega,\nu)$. We establish the converse when $\omega$ has a $p$-small highest weight.

\medskip

\noindent {\it{A local mod p Langlands correspondence?}} It is widely believed that there should be a local mod $p$ correspondence, which is compatible with $\frak{B}$ above.
At the very least, there should be map $\bar{\rho}\mapsto b(\bar{\rho})$ (on equivalence classes)
$$
\left\{ \begin{matrix} \text{continuous $\bar{\rho}:\Gamma_K \rightarrow \GL_n(\F_L)$}  \end{matrix} \right\}
\overset{?}{\longrightarrow}
$$
$$
\left\{ \begin{matrix} \text{$\F_L$-spaces $b$ endowed with} \\ \text{admissible $\GL_n(K)$-action}  \end{matrix} \right\},
$$
which takes an irreducible $\bar{\rho}$ to an irreducible $b(\bar{\rho})$. One should be able to recover $\bar{\rho}$ from $b(\bar{\rho})$. Furthermore, for any lift 
$\rho: \Gamma_K \rightarrow \GL_n(\mathcal{O}_L)$ of $\bar{\rho}$, there ought to be at least a canonical nonzero map (perhaps an embedding),
$$
\frak{B}(\rho\otimes L)^{\circ}\otimes \F_L \overset{?}{\longrightarrow} b(\bar{\rho}).
$$
At the end, we define an admissible mod $p$ representation $b(\bar{\rho}_x)$ of $G(\Q_p)$, when $\bar{\rho}_x$ is irreducible, which we like to think of as a candidate for $\otimes_{v|p}b(\bar{\rho}_x|_{\Gamma_{\K_{\tilde{v}}}})$, at least when these local restrictions above $p$ are irreducible. However, we cannot say much about it. Even showing that $b(\bar{\rho}_x)\neq 0$ appears to be a difficult problem, closely tied to the conjectural Ihara lemma for unitary groups (Conjecture B in [CHT]). We formulate a precise local-global compatibility conjecture for $b(\bar{\rho}_x)$, and explain how it {\it{implies}} Ihara's lemma. This makes heavy use of [EH].

\medskip

\noindent {\it{Acknowledgements}}. With great pleasure and admiration, I want to acknowledge the impact of the visions of Breuil and Emerton on this work. I would like to thank 
both of them heartfully for discretely pointing out some (very) embarrassing misconceptions of mine, the first time I spoke about it (in hindsight prematurely).
Also, conversations and correspondence with Chojecki, Dospinescu, Helm, Herzig, Newton, Ramakrishnan, and Shin, have been a great help and a source of inspiration. I am grateful to Minguez for clarifying Corollary 4.2 of [Min] to me, and for pointing to the reference [CHL]. Thanks are due to Caraiani, for making a preliminary version of [Car] available to me before circulating it.

\section{Automorphic Galois representations}

We start out by summarizing what is currently known about attaching Galois representations to automorphic representations of definite unitary groups. Due to the work of many people, we now have an almost complete understanding of this, and below we merely navigate the existing literature. We claim no originality in this section. Our goal is simply to state the precise result. Particularly, we want to emphasize the local-global compatibility at $p=\ell$, recently proved in [BGGT] and [Car], which is fundamental for this paper.

\subsection{Definite unitary groups}

Throughout this article, we fix a totally real field $F$, and a CM extension $\K$. We let $c$ denote the non-trivial element of $\Gal(\K|F)$. The places of $F$ will usually be denoted by $v$, those of $\K$ by $w$. We are interested in outer forms $U$ of $\GL(n)_{F}$, which become an inner form $D^{\times}$ over $\K$. Here $D$ is a central simple $\K$-algebra, of $\dim_{\K}(D)=n^2$. These forms are unitary groups $U=U(D,\star)$, where $\star$ is an anti-involution on $D$ of the second kind ($\star|_{\K}=c$). Thus, for 
any $F$-algebra $R$,
$$
U(R)=\{x \in (D\otimes_F R)^{\times}: xx^{\star}=1\}.
$$
We will always assume from now on that $U(F\otimes_{\Q}\R)$ is {\it{compact}}. Thus, by making a choice of a CM-type $\Phi$, the group may be identified with
$U(n)^{\Hom(F,\R)}$ (up to conjugation). It will be convenient to work over the rationals, and introduce $G=\text{Res}_{F|\Q}(U)$. With the same $\Phi$ one identifies $G(\C)$ with
$\GL_n(\C)^{\Hom(F,\R)}$.

\subsection{Weights of automorphic representations}

Following standard notation in the subject, $(\Z^n)_+^{\Hom(\K,\C)}$ will denote the set of tuples $a=(a_{\tau})_{\tau\in \Hom(\K,\C)}$, where each $a_{\tau}=(a_{\tau,j})$ itself is a decreasing tuple,
$$
a_{\tau}=(a_{\tau,1}\geq a_{\tau,2}\geq \cdots \geq a_{\tau,n}),
$$
of integers. In the obvious way, we can identify $a_{\tau}$ with a dominant weight for $\GL(n)$, relative to the upper triangular Borel. We say $a_{\tau}$ is {\it{regular}} if all the inequalities above are strict. We say $a$ is regular if $a_{\tau}$ is regular for {\it{some}} $\tau$.

\medskip

\noindent Now, let $\pi=\pi_{\infty}\otimes \pi_f$ be an automorphic representation of $U(\A_F)$. We will define what it means for $\pi$ to have weight $a$: Every embedding $\tau: \K \hookrightarrow \C$ restricts to a $\sigma: F \hookrightarrow \R$, which corresponds to an infinite place $v=v(\sigma)$ of $F$. With this notation, $\tau$ identifies $U(F_v)\simeq U(n)$, under which $\pi_v$ should be equivalent to the contragredient $\breve{V}_{a_{\tau}}$, or rather its restriction. Here $V_{a_{\tau}}$ is the irreducible algebraic representation of $\GL_n(\C)$ of highest weight $a_{\tau}$. 

\medskip

\noindent {\it{Remark}}. We must have $V_{a_{\tau c}}=\breve{V}_{a_{\tau}}$. In other words, $a_{\tau c,j}=-a_{\tau,n+1-j}$.

\subsection{Associating Galois representations}

We have introduced enough notation, in order to formulate the following main result, the foundation for our work. As mentioned already, this is the culmination of collaborative efforts of a huge group of outstanding mathematicians, as will become clear below.

\begin{thm}
Choose a prime $p$, and an isomorphism $\iota: \C \overset{\sim}{\longrightarrow} \bar{\Q}_p$.
Let $\pi$ be an automorphic representation of $U(\A_F)$ such that $\pi_{\infty}$ has {\bf{regular}} weight $a$. Then there exists a unique continuous semisimple Galois representation
$$
\rho_{\pi,\iota}: \Gamma_{\K}=\Gal(\bar{\Q}|\K)\rightarrow \GL_n(\bar{\Q}_p)
$$
such that the following properties are satisfied:
\begin{itemize}
\item[(a)] $\breve{\rho}_{\pi,\iota}\simeq \rho_{\pi,\iota}^c \otimes \epsilon_{cyc}^{n-1}$,
\item[(b)] For {\bf{every}} finite place $v$, and every $w|v$ (even those above $p$),
$$
WD(\rho_{\pi,\iota}|_{\Gamma_{\K_w}})^{F-ss}\simeq \iota \text{rec}(BC_{w|v}(\pi_v)\otimes |\det|_w^{(1-n)/2}),
$$
whenever $BC_{w|v}(\pi_v)$ is defined: If $\pi_v$ is unramified or $v=ww^c$ splits.
\item[(c)] $\rho_{\pi,\iota}|_{\Gamma_{\K_w}}$ is potentially semistable for all $w|p$, with Hodge-Tate numbers
$$
HT_{\iota\tau}(\rho_{\pi,\iota}|_{\Gamma_{\K_w}})=\{a_{\tau,j}+(n-j): j=1,\ldots,n\}
$$
for every $\tau:\K \hookrightarrow \C$ such that $\iota\tau$ lies above $w$. A word about our normalization here;
$\rho_{\pi,\iota}\otimes_{\iota\tau,\K_w}\C_{\K_w}(i)$ has no $\Gamma_{\K_w}$-invariants unless $i$ is of the above form, in which case they form a line.
Thus $HT_{\iota\tau}(\epsilon_{cyc})=\{-1\}$.
\end{itemize}
\end{thm}

\noindent {\it{Proof}}. Ngo's proof of the fundamental lemma makes functoriality available in a slew of cases. In particular, weak base change from any unitary group associated with $\K|F$ to $\GL_n(\A_{\K})$ has matured. Building on work of Clozel and Labesse, White was recently able to work out the cohomological case completely [Whit]. In our given setup, 
$\pi_v$ is automatically discrete series for all $v|\infty$, in which case Theorem 6.1 in [Whit], or rather the pertaining remarks 6.2 and 6.3, yields an automorphic representation
$$
\Pi=\Pi_1\boxplus \cdots \boxplus \Pi_t
$$
of $\GL_n(\A_{\K})$, which is an isobaric sum of mutually non-isomorphic conjugate self-dual cuspidal automorphic representations $\Pi_i$ of some $\GL_{n_i}(\A_{\K})$ such that
$$
\Pi_w=\text{BC}_{w|v}(\pi_v)
$$
for all $w|v$, where $v$ is split or archimedean, or $\pi_v$ is unramified. The {\it{regularity}} of $\pi_{\infty}$ ensures that the $\Pi_i$ are cuspidal (as opposed to just discrete), which in turn implies the previous equality at the archimedean places $w|v$. Let us spell it out in that case: Fix an embedding $\tau:\K \hookrightarrow \C$ inducing $\K_w \simeq \C$. Then,
$$
\phi_{\Pi_w}: \C^*=W_{\C}\simeq W_{\K_w}\rightarrow \GL_n(\C)
$$
maps
$$
z \mapsto \begin{pmatrix}(z/\bar{z})^{-h_1+\frac{n-1}{2}} & & \\ & \ddots & \\ & & (z/\bar{z})^{-h_n+\frac{n-1}{2}}\end{pmatrix},
$$
for certain $h_j \in \Z$, which are given in terms of the weight by $h_j=a_{\tau,j}+(n-j)$. This last formula is worked out in [BeCl] for example. These $h_j$ are distinct, so each $\Pi_i\otimes|\det|_{\K}^{(n_i-n)/2}$ is regular algebraic, {\it{essentially}} conjugate self-dual, and cuspidal. By Theorem A of [BGGT], and the references therein, we can associate a Galois representation $r_{\Pi_i,\iota}$ to it, satisfying the analogous properties of (a) to (c). As a remark, in loc. cit. local-global compatibility {\it{at}} $p$ is proved assuming Shin-regularity, which is much weaker than regularity. In any case, the Shin-regularity assumption is currently being removed by Caraiani as part of her 2012 Harvard Ph.D thesis. It is then straightforward to check that the representation
$$
\rho_{\pi,\iota}=r_{\Pi_1,\iota}\oplus \cdots \oplus r_{\Pi_t,\iota}
$$
has the desired properties. It is uniquely determined by (b) by Tchebotarev. $\square$

\medskip

\noindent {\it{Remark}}. It appears within reach to extend the previous argument to the {\it{irregular}} case. By [Whit], one still has a weak base change $\boxplus_{i=1}^t\Pi_i$, but the $\Pi_i$ are only {\it{discrete}}, not cuspidal. By Shapiro's lemma in $(\frak{g},K)$-cohomology, these $\Pi_i$ should still be cohomological (of Speh type). By the Moeglin-Waldspuger description of the discrete spectrum of $\GL(n_i)$, one can in turn express each $\Pi_i$ as an isobaric sum of cusp forms, with which one can associate Galois representations. After having consulted several experts in the field, we are quite optimistic about this line of argument, and that the ubiquitous regularity assumption (appearing throughout this paper) can safely be dropped. However, we have not made any serious attempt to work out the details. We are hopeful that forthcoming joint work of Kaletha-Shin-White should provide the strengthenings of [Whit] needed.

\subsection{The Breuil-Schneider recipe}

As described in detail in the introduction, to a potentially semistable representation $\rho:\Gamma_{\K_w}\rightarrow \GL_n(\bar{\Q}_p)$, with distinct Hodge-Tate numbers, Breuil and Schneider attach a locally algebraic representation $\text{BS}(\rho)$ of $\GL_n(\K_w)$ on a $\bar{\Q}_p$-vector space. The algebraic part is given by all the Hodge-Tate numbers $\text{HT}(\rho)$ (for varying embeddings $\K_w \hookrightarrow \bar{\Q}_p$), and the smooth part is given by the Weil-Deligne representation $\text{WD}(\rho)$, or rather its Frobenius-semisimplification, via the classical local Langlands correspondence, slightly modified in the non-generic case. They conjecture the mere existence of an invariant norm on $\text{BS}(\rho)$. In fact,
Conjecture 4.3 in [BrSc] also predicts a converse, which has been proved in [Hu]. Our goal is to first prove the Breuil-Schneider conjecture for $\rho=\rho_{\pi,\iota}|_{\Gamma_{\K_w}}$, for any place $w$ of $\K$ above $p$. We will achieve this below. For now, we will compute  $\text{BS}(\rho)$ explicitly in this situation.

\medskip

\noindent In fact, we prefer to use a slightly different normalization: There is a choice of a {\it{sign}} involved in the recipe on p. 16 in [BrSc]. Instead of twisting by $|\det|_w^{(1-n)/2}$, we prefer to twist by  $|\det|_w^{(n-1)/2}$ to make it more compatible with the previous notation. Consequently, $\text{BS}(\rho)$ becomes twisted by an integral character.

\begin{df}
For $a \in \Q_p^*$ we let $a^{\times}=a|a|_p$ denote its unit part. We introduce
$$
\mu: \GL_n(\K_w)\overset{\det}{\longrightarrow} \K_w^*\overset{N_{\K_w|\Q_p}}{\longrightarrow}\Q_p^*\longrightarrow \Q_p^*/p^{\Z}\simeq \Z_p^{\times}.
$$
That is, $\mu(g)=N_{\K_w|\Q_p}(\det g)^{\times}$. We will normalize $\text{BS}(\rho)$ as follows:
$$
\widetilde{\text{BS}}(\rho):=\text{BS}(\rho)\otimes_{\bar{\Q}_p}\mu^{n-1}.
$$
(Of course, this has an invariant norm precisely when $\text{BS}(\rho)$ does.)
\end{df}

\begin{lem}
$\widetilde{\text{BS}}(\rho)=\text{BS}(\rho(n-1))$.
\end{lem}

\noindent {\it{Proof}}. Note that the character $a \mapsto a^{\times}$ (which maps $p \mapsto 1$, and is the identity on $\Z_p^{\times}$) corresponds to the $p$-adic cyclotomic character $\chi_{cyc}: \Gamma_{\Q_p}\rightarrow \Z_p^{\times}$ via local class field theory $\Q_p^*\rightarrow\Gamma_{\Q_p}^{ab}$. For any $p$-adic field $K$, it follows that
$BS(\chi_{cyc})$ is simply the character $a \mapsto N_{K|\Q_p}(a)^{\times}$. Consequently, $\widetilde{\text{BS}}(\rho)=\text{BS}(\rho\otimes \chi_{cyc}^{n-1})$. $\square$

\medskip

\noindent We compute it in the automorphic case. Given the local-global compatibility results of [BGGT] (generalized in [Car]), this is basically just "bookkeeping".

\begin{cor}
Let $\pi$ be an automorphic representation of $U(\A_F)$ of regular weight $a$. Assume $\rho_{\pi,\iota}$ is absolutely {\bf{irreducible}} (as a representation of the full Galois group $\Gamma_{\K}$). Let $v|p$ be a place of $F$, either split in $\K$, or such that $\pi_v$ is unramified. Then, for any place $w|v$ of $\K$, we have
$$
\widetilde{\text{BS}}(\rho_{\pi,\iota}|_{\Gamma_{\K_w}})=(\otimes_{\sigma: \K_w \hookrightarrow \bar{\Q}_p} \breve{V}_{a_{\iota^{-1}\sigma}}) 
\otimes_{\bar{\Q}_p}(BC_{w|v}(\pi_v)\otimes_{\C,\iota}\bar{\Q}_p).
$$
(We abuse notation a bit and let $V_{a_{\tau}}$ denote the irreducible algebraic representation of $\GL_n(\bar{\Q}_p)$ of highest weight $a_{\tau}$; as opposed to the {\it{complex}} representation from earlier chapters)
\end{cor}

\noindent {\it{Proof}}.  What is denoted $\pi^{\text{unit}}$ in [BrSc] equals $BC_{w|v}(\pi_v)\otimes |\det|_w^{(1-n)/2}$ in our case (more precisely, $\otimes_{\C,\iota}\bar{\Q}_p$). 
When it is generic, the smooth part of $\text{BS}(\rho)$ is
$$
\pi^{\text{unit}}\otimes_{\bar{\Q}_p} |\det|_w^{(1-n)/2}=(BC_{w|v}(\pi_v)\otimes |\det|_w^{1-n})\otimes_{\C,\iota}\bar{\Q}_p.
$$
In the non-generic case, $\pi^{\text{unit}}$ has to be replaced by a certain parabolically induced representation. However, if we assume $\rho_{\pi,\iota}$ is (globally) irreducible,
we see that $\Pi=\text{BC}_{\K|F}(\pi)$ must be cuspidal, and in particular $\Pi_w$ is generic. The algebraic part of $\text{BS}(\rho)$ is constructed out of the Hodge-Tate numbers:
What is denoted $i_{j,\sigma}$ in [BrSc], for an embedding $\sigma:\K_w \hookrightarrow \bar{\Q}_p$, equals $a_{\tau,n+1-j}+(j-1)$ in our notation, where
$\sigma=\iota\tau$. In (8) on p. 17 of [BrSc], the numbers become
$$
b_{\tau,j}:=-i_{n+1-j,\sigma}-(j-1)=-a_{\tau,j}-(n-1).
$$
Breuil-Schneider's $\rho_{\sigma}$ is the irreducible algebraic representation of $\GL_n(\bar{\Q}_p)$ of highest weight $b_{\tau,1}\leq \cdots\leq b_{\tau,n}$ relative to the {\it{lower}} triangular Borel. Relative to the upper triangular Borel, $\rho_{\sigma}$ has highest weight $b_{\tau,n}\geq \cdots\geq b_{\tau,1}$, so that
$\rho_{\sigma}\simeq \breve{V}_{a_{\tau}}\otimes \det^{1-n}$ (more precisely, $\otimes_{\C,\iota}\bar{\Q}_p$). Altogether, the algebraic part is
$$
\xi=\otimes_{\sigma} \rho_{\sigma}\simeq \otimes_{\tau|w} (\breve{V}_{a_{\tau}}\otimes {\det}^{1-n})
$$
(the tensor product ranging over $\tau:\K \hookrightarrow \C$ such that $\iota\tau$ induces $w$). Here we abuse notation a bit, and use $V_{a_{\tau}}$ to denote the irreducible algebraic representation of $\GL_n(\bar{\Q}_p)$ of highest weight $a_{\tau}$. As a representation of $\GL_n(\K_w)$, embedded diagonally in 
$\prod_{\sigma: \K_w \hookrightarrow \bar{\Q}_p}\GL_n(\bar{\Q}_p)$, the algebraic part becomes
$$
\xi=(\otimes_{\sigma: \K_w \hookrightarrow \bar{\Q}_p} \breve{V}_{a_{\iota^{-1}\sigma}}) \otimes (N_{\K_w|\Q_p}\circ \det)^{1-n},
$$
which yields the result. $\square$

\section{Completed Cohomology}

In this section we will prove the Breuil-Schneider conjecture, 4.3 in [BrSc], for the potentially semistable representations $\rho=\rho_{\pi,\iota}|_{\Gamma_{\K_w}}$ above. This will make heavy use of ideas of Emerton introduced in [Emer]. The basic idea is to view $\widetilde{\text{BS}}(\rho)$ as a component of the $p$-adic automorphic representation $\tilde{\pi}=\tilde{\pi}_p\otimes \pi_f^p$ attached to $\pi$, which in turn embeds into the completed cohomology $\tilde{H}^0$ for $G$. 

\subsection{The $p$-adic automorphic representation}

We keep our automorphic representation $\pi$ of $U(\A_F)$ of regular weight $a$. Recall that we introduced the group $G=\text{Res}_{F|\Q}(U)$.
Interchangeably, below we will view $\pi$ as an automorphic representation of $G(\A)$. We will follow p. 52 in [Emer] in attaching a $p$-adic automorphic representation to $\pi$.
(The $\G$ there will be our $G$, and $F$ there will be $\Q$.) This can be done for $W$-allowable $\pi$, where $W$ is an irreducible algebraic representation of $G(\C)$, which in this case (where $G$ is compact at infinity) simply means $\pi_{\infty}\simeq W|_{G(\R)}$. See Definition 3.1.3 in [Emer]. 

\medskip

\noindent To make this more explicit, in terms of the weight $a$, we need to make some identifications. Let us choose a CM-type $\Phi$. For each $\sigma:F \hookrightarrow \R$
we let $\tilde{\sigma}$ denote its lift in $\Phi$. Thus the two extensions to $\K$ are $\{\tilde{\sigma},\tilde{\sigma}^c\}$. Via the choice of $\Phi$,
$$
\text{$G(\C)\overset{\sim}{\longrightarrow}_{\Phi} \GL_n(\C)^{\Hom(F,\R)}$, $\y$ $G(\R)\overset{\sim}{\longrightarrow}_{\Phi} U(n)^{\Hom(F,\R)}$.}
$$
We immediately infer that $W\simeq \otimes_{\sigma \in \Hom(F,\R)} \breve{V}_{a_{\tilde{\sigma}}}$ under these identifications. Via 
$\iota: \C \overset{\sim}{\longrightarrow} \bar{\Q}_p$ we identify $W$ with an algebraic representation of $G(\bar{\Q}_p)$.
$$
G(\bar{\Q}_p)\overset{\sim}{\longrightarrow}_{\Phi} \prod_{v|p}\GL_n(\bar{\Q}_p)^{\Hom(F_v,\bar{\Q}_p)}
$$
allows us to factor our $p$-adic $W$ accordingly, as $W\simeq \otimes_{v|p}W_v$, where we let
$$
W_v=\otimes_{\sigma \in \Hom(F,\R): \sigma|v}\breve{V}_{a_{\tilde{\sigma}}}.
$$
In the same vein, $G(\Q_p)=\prod_{v|p}U(F_v)$. To go any further, from this point on we assume every $v|p$ splits in $\K$, and that $D_w\simeq M_n(\K_w)$ for each divisor $w|v$.
Then $U(F_v)\overset{\sim}{\longrightarrow} \GL_n(\K_w)$, defined up to conjugation. If we assume (as we may) that our CM-type $\Phi$ is ordinary at $\iota$, in the sense of [Katz], $\Phi$ singles out a place $\tilde{v}$ of $\K$ above each $v|p$ of $F$. With this selection of places at hand,
$$
G(\Q_p) \overset{\sim}{\longrightarrow} \prod_{v|p} \GL_n(\K_{\tilde{v}}).
$$
Moreover, the inclusion into $G(\bar{\Q}_p)$ corresponds to the diagonal embeddings,
$$
\GL_n(\K_{\tilde{v}})=\GL_n(F_v)\hookrightarrow \GL_n(\bar{\Q}_p)^{\Hom(F_v,\bar{\Q}_p)}.
$$
The following is Definition 3.1.5 in [Emer], except that we are working with representations over $\bar{\Q}_p$ instead of descending to a finite extension of $\Q_p$.

\begin{df}
The classical $p$-adic automorphic representation of $G(\A_f)$ over $\bar{\Q}_p$ attached to the $W$-allowable automorphic representation $\pi$ of $G(\A)$ is
$$
\text{$\tilde{\pi}:=\tilde{\pi}_p\otimes_{\bar{\Q}_p}\pi_f^p$, $\y$ $\tilde{\pi}_p:=W \otimes_{\bar{\Q}_p}\pi_p$.}
$$ 
Here $G(\Q_p)$ acts diagonally on $W \otimes_{\bar{\Q}_p}\pi_p$, and $G(\A_f^p)$ acts through the second factor $\pi_f^p$. (Abusing notation, we write $\pi_p$ instead of $\pi_p\otimes_{\C,\iota}\bar{\Q}_p$ and so on.)
\end{df}

\noindent At each $v|p$ we introduce $\tilde{\pi}_v=W_v \otimes_{\bar{\Q}_p}\text{BC}_{\tilde{v}|v}(\pi_v)$, a locally algebraic representation of $\GL_n(\K_{\tilde{v}})$, which depends on the choice of an ordinary CM-type $\Phi$. Moreover, $\tilde{\pi}_p\simeq \otimes_{v|p}\tilde{\pi}_v$ under the isomorphism $G(\Q_p)\simeq \prod_{v|p}\GL_n(\K_{\tilde{v}})$.

\medskip

\noindent This leads to the main result of this section.

\begin{prop}
Hypothesis: Every $v|p$ of $F$ splits in $\K$, and $D_w\simeq M_n(\K_w)$ for all $w|v$. For each $v|p$ of $F$ pick a place $\tilde{v}|v$ of $\K$ (this amounts to choosing an $\iota$-ordinary CM-type). Let $\pi$ be an automorphic representation of $U(\A_F)$ of regular weight, and assume
$\rho_{\pi,\iota}$ is (globally) irreducible. Then, for all $v|p$ of $F$,
$$
\widetilde{\text{BS}}(\rho_{\pi,\iota}|_{\Gamma_{\K_{\tilde{v}}}})\simeq \tilde{\pi}_v,
$$
which embeds into $\tilde{\pi}|_{\GL_n(\K_{\tilde{v}})}$ (where we restrict via $U(F_v)\overset{\sim}{\longrightarrow}\GL_n(\K_{\tilde{v}})$).
\end{prop}

\noindent {\it{Proof}}. This follows from the preceding discussion, combined with the computation of the Breuil-Schneider representation in Corollary 1 above. $\square$

\subsection{Algebraic modular forms}

We will study the space of modular forms for $G$ of a given weight. To put things in a broader perspective, we will use the cohomological framework of [Emer], although we will only work with $H^0$, which is explicit, and of a combinatorial nature. In our situation,
$G(\R)$ is compact and connected, so things simplify tremendously, and we only have cohomology in degree zero. Indeed, for every compact open subgroup $K\subset G(\A_f)$, 
the corresponding arithmetic quotient is a {\it{finite}} set:
$$
Y(K)=G(\Q)\backslash G(\A_f)/K.
$$
An irreducible algebraic representation $W$ of $G(\C)$ defines a local system $\mathcal{V}_W$ on each $Y(K)$, and $H^0(Y(K),\mathcal{V}_{\breve{W}})$ is identified with the space of modular forms, of level $K$, and weight $W$. That is, all functions $f: G(\A_f)\rightarrow \breve{W}$, which are $K$-invariant on the right, such that $f(\gamma g)=\gamma f(g)$ for all elements $\gamma \in G(\Q)$.
$$
H^0(\mathcal{V}_{\breve{W}}):=\underset{K}{\varinjlim}  H^0(Y(K),\mathcal{V}_{\breve{W}})\simeq \oplus_{\pi: \pi_{\infty}\simeq W}m_G(\pi)\pi_f
$$
is then a smooth admissible semisimple representation of $G(\A_f)$, which we wish to suitably $p$-adically complete. Via our choice of $\iota: \C \overset{\sim}{\longrightarrow} \bar{\Q}_p$, we will view $W$ as a representation of $G(\bar{\Q}_p)$ and so on. Occasionally it will be convenient to work over a field $E\subset \bar{\Q}_p$, finite over $\Q_p$. 
It suffices to take $E$ large enough, so that it contains the image of every embedding $F \hookrightarrow \bar{\Q}_p$. In that case $G$ splits over $E$, and by highest weight theory
$W$ may be defined over $E$. Thus, from now on, $H^0(\mathcal{V}_{\breve{W}})$ is an $E$-vector space with a smooth admissible $G(\A_f)$-action.

\begin{df}
For each tame level $K^p\subset G(\A_f^p)$, following [Emer], we introduce
$$
H^0(K^p,\mathcal{O}_E/\varpi_E^s):=\underset{K_p}{\varinjlim} H^0(Y(K_pK^p),\mathcal{O}_E/\varpi_E^s),
$$
and
$$
\tilde{H}^0(K^p):=E \otimes_{\mathcal{O}_E}\underset{s}{\varprojlim} H^0(K^p,\mathcal{O}_E/\varpi_E^s).
$$
The latter is an $E$-Banach space with a unitary $G(\Q_p)$-action, commuting with the action of the Hecke algebra $\mathcal{H}(K^p)$ of compactly supported $K^p$-biinvariant $E$-valued functions
on $G(\A_f^p)$. In fact, it becomes a Banach module over the completion $\hat{\mathcal{H}}(K^p)$. Also,
$$
\tilde{H}^0:=\underset{K^p}{\varinjlim} \tilde{H}^0(K^p),
$$
a locally convex $E$-vector space with an action of $G(\A_f)$.
\end{df}

\noindent In our simple setup, they can all be realized very explicitly. For example,
$$
\text{$\tilde{H}^0(K^p)=\{\text{continuous $Y(K^p)\overset{f}{\rightarrow} E$}\}$, $\y$ $Y(K^p)=\underset{K_p}{\varprojlim} Y(K_pK^p)$,}
$$
with the supremum-norm $\|\cdot\|$. The {\it{key}} ingredient we will use is the isomorphism:

\begin{lem}
For any absolutely irreducible algebraic representation $W$, 
$$
W \otimes_E H^0(\mathcal{V}_{\breve{W}}) \overset{\sim}{\longrightarrow} (\tilde{H}^0)_{\text{$W$-alg}}.
$$ 
(When $W$ is only irreducible over $E$, tensor over $\End_{\frak{g}}(W)$, where $\frak{g}=\text{Lie}G(\Q_p)$.)
\end{lem}

\noindent {\it{Proof}}. This is Corollary 2.2.25 in [Emer] (also spelled out in [Sor] for $H^0$). Let us briefly sketch the main idea. For any tame level $K^p$, one shows that
$$
W \otimes_E H^0(K^p,\mathcal{V}_{\breve{W}})=\underset{K_p}{\varinjlim} W \otimes_E H^0(Y(K_pK^p),\mathcal{V}_{\breve{W}})) \overset{\sim}{\longrightarrow} \tilde{H}^0(K^p)_{\text{$W$-alg}}.
$$ 
This goes as follows: $H^0(Y(K_pK^p),\mathcal{V}_{\breve{W}}))$ is a space of classical $p$-adic modular forms, and it is an easy exercise to identify it with
$\Hom_{K_p}(W, \tilde{H}^0(K^p))$. Now,
$$
W \otimes_E \Hom_{K_p}(W, \tilde{H}^0(K^p))\overset{eval.}{\longrightarrow} \tilde{H}^0(K^p)
$$
is injective since $W$ is absolutely irreducible, even when restricted to $K_p$ (which is Zariski dense). The image of this evaluation map is the $W$-isotypic subspace of $\tilde{H}^0(K^p)$. As $K_p$ varies, the maps are compatible, and produces a map out of the direct limit onto $\tilde{H}^0(K^p)_{\text{$W$-alg}}$ as desired. $\square$

\medskip

\noindent {\it{Remark}}. For higher degree cohomology $H^i$ there is an analogous canonical $G(\A_f)$-equivariant map, which occurs as the edge map of a certain spectral sequence, but it is not known to be injective for groups other than $\GL(2)_{\Q}$ (and groups $G$ which are compact at infinity mod center). Injectivity is what makes the whole machinery of [Emer] work, see his Theorem 0.7 and Proposition 2.3.8 on p. 47, for example. In particular, it is available in our case, where $G(\R)$ is compact. In general, one would have to localize the spectral sequence at a "cohomologically" non-Eisenstein maximal ideal $\m$ (which means it does not contribute to mod $p$ cohomology outside the middle degree). This is expected to hold when the Galois representation $\bar{\rho}_{\m}$ is absolutely irreducible, but this is difficult to show. Partial results are now available for $U(2,1)$. 

\medskip

\noindent From the previous discussion,  we get decompositions of completed cohomology:

\begin{prop}
\begin{itemize}
\item[(1)] $\bar{\Q}_p\otimes_E (\tilde{H}^0)_{\text{$W$-alg}} \simeq \oplus_{\pi:\pi_{\infty}\simeq W} m_G(\pi) \tilde{\pi}$,
\item[(2)] $\bar{\Q}_p\otimes_E \tilde{H}^0(K^p)_{\text{$W$-alg}} \simeq \oplus_{\pi:\pi_{\infty}\simeq W} m_G(\pi) (\tilde{\pi}_p\otimes_{\bar{\Q}_p} (\pi_f^p)^{K^p})$.
\end{itemize}
\end{prop}

\noindent Now, suppose $\frak{h}\subset \mathcal{H}(K^p)$ is a central subalgebra. It then acts on $(\pi_f^p)^{K^p}$ by a character $\lambda_{\pi}: \frak{h}\rightarrow \bar{\Q}_p$.
Conversely, say we start out with $\lambda: \frak{h}\rightarrow \bar{\Q}_p$. Then,
$$
\bar{\Q}_p\otimes_E \tilde{H}^0(K^p)_{\text{$W$-alg}}^{\frak{h}=\lambda} \simeq \oplus_{\pi:\pi_{\infty}\simeq W, \lambda_{\pi}=\lambda} m_G(\pi) (\tilde{\pi}_p\otimes_{\bar{\Q}_p} (\pi_f^p)^{K^p}).
$$
As always, we assume $W$ has {\it{regular}} weight, so we know how to attach Galois representations. If $\frak{h}$ contains the spherical part $\mathcal{H}(K^p)^{\text{sph}}$, all the $\pi$ contributing to the right-hand side have the same Galois representation $\rho_{\lambda}$, by Chebotarev, which we assume is {\it{irreducible}}. By Proposition 1, we may factor the above,
$$
\bar{\Q}_p\otimes_E \tilde{H}^0(K^p)_{\text{$W$-alg}}^{\frak{h}=\lambda} \simeq (\otimes_{v|p} \widetilde{BS}(\rho_{\lambda}|_{\Gamma_{\K_{\tilde{v}}}})) \otimes_{\bar{\Q}_p}
(\oplus_{\pi:\pi_{\infty}\simeq W,\lambda_{\pi}=\lambda} m_G(\pi)\pi_f^p)^{K^p}.
$$
This has the form of a $G(\Q_p)\simeq \prod_{v|p}\GL_n(\K_{\tilde{v}})$-representation tensor an $\mathcal{H}(K^p)$-module. In particular, since $\tilde{H}^0(K^p)$ carries
a $G(\Q_p)$-invariant norm, we finally deduce the Breuil-Schneider conjecture for automorphic Galois representations:

\begin{thm}
If $\pi$ is an automorphic representation of $U(\A_F)$, of regular weight, such that $\rho_{\pi,\iota}$ is irreducible. Then $\text{BS}(\rho_{\pi,\iota}|_{\Gamma_{\K_w}})$ admits
a $\GL_n(\K_w)$-invariant norm, for all places $w|p$ of $\K$.
\end{thm}

\noindent The discussion leading up to the Theorem strongly suggests a better formulation in terms of eigenvarieties. We will employ this machinery in the next Chapter.

\section{Eigenvarieties}

Eigenvarieties are rigid analytic spaces interpolating Hecke eigensystems occurring in spaces of automorphic forms, of varying weight. Historically, the first example is the Coleman-Mazur eigencurve for $\GL(2)_{\Q}$, revisited by Buzzard, Emerton, Urban, and others. There are different constructions for any reductive group $G$, which each have their drawbacks and limitations. When $G(\R)$ is compact, however, the theory is in good shape, and all constructions are compatible. Below we will combine the approach of [Emer] with that of [Chen] (for arbitrary totally real $F$) extending parts of  [BeCh] (when $F=\Q$).

\subsection{The classical points}

By our standing hypotheses, $G_{\Q_p}\simeq \prod_{v|p}\text{Res}_{\K_{\tilde{v}}|\Q_p}\GL(n)$ is quasi-split, and we pick the Borel pair $(B,T)$, defined over $\Q_p$, corresponding to the product of the upper triangular pairs in each $\GL_n(\K_{\tilde{v}})$. 

\medskip

\noindent As in [Emer], let $\hat{T}$ denote the {\it{weight}} space. That is, the rigid analytic variety over $E$ (a subfield of $\bar{\Q}_p$, finite over $\Q_p$, but large enough so that $G$ splits over $E$) which parametrizes the locally analytic characters on $T(\Q_p)$. In other words, 
$$
\hat{T}(A)=\Hom_{la}(T(\Q_p),A^{\times})
$$
for any affinoid $E$-algebra $A$. It comes with a universal map $T(\Q_p)\rightarrow \mathcal{O}(\hat{T})^{\times}$.

\medskip

\noindent The eigenvariety depends on the choice of {\it{tame}} level $K^p\subset G(\A_f^p)$, which we will always assume is decomposable as $\prod_{v\nmid p}K_v$, where $K_v$ is a compact open subgroup of $U(F_v)$, which is hyperspecial for all but finitely many $v$. Say, for all $v \notin S(K^p)$. Correspondingly, the Hecke algebra factors as a tensor product,
$$
\mathcal{H}(K^p)=\otimes_{v\nmid p} \mathcal{H}(K_v)=\mathcal{H}(K^p)^{\text{ram}}\otimes_E \mathcal{H}(K^p)^{\text{sph}}.
$$
Here $\mathcal{H}(K^p)^{\text{sph}}=\otimes_{v\notin S(K^p)} \mathcal{H}(K_v)$ sits as a central subalgebra of $\mathcal{H}(K^p)$, hence it acts by a character on $\pi_f^{K^p}$, for any automorphic $\pi$ with $K^p$-invariants.

\medskip

\noindent We now make precise which points we wish to interpolate by an eigenvariety.

\begin{df}
Let $E(0,K^p)_{cl}\subset (\hat{T}\times \text{Spec}  \mathcal{H}(K^p)^{\text{sph}})(\bar{\Q}_p)$ be the subset of pairs $x=(\chi,\lambda)$ for which there exists an irreducible 
$G(\A_f)$-subquotient $\pi_f$ of $\bar{\Q}_p\otimes_E H^0(\mathcal{V}_{\breve{W}})$, where $W$ is an irreducible algebraic representation of $G_E$, such that
\begin{itemize}
\item[(a)] $\chi=\psi\theta$, where $\psi$ is the highest weight of $W$ (relative to $B$), and $\theta$ is a smooth character of $T(\Q_p)$ such that $\pi_p \hookrightarrow \text{Ind}_{B(\Q_p)}^{G(\Q_p)}(\theta)$,
\item[(b)] $\pi_f^{K^p}\neq 0$, and $\mathcal{H}(K^p)^{\text{sph}}$ acts on it via $\lambda$.
\end{itemize}
\end{df}

\noindent This is the definition, and notation, used on p. 5 in [Emer].

\subsection{Eigenvariety conventions}

Emerton defines the degree zero cohomological eigenvariety of $G$, of tame level $K^p$, to be the rigid analytic closure of  $E(0,K^p)_{cl}$ in
$\hat{T}\times \text{Spec}  \mathcal{H}(K^p)^{\text{sph}}$. By the uniqueness part of Theorem 1.6 in [Chen], it coincides with the eigenvariety defined there. We will intertwine the two points of view. Thus, with $E(0,K^p)_{cl}$ is associated a quadruple $(\X,\chi,\lambda,X_{cl})$, consisting of the following data:

\begin{itemize}
\item $\X_{/E}$ is an equi-dimensional reduced rigid analytic variety,
\item $\chi: \X \rightarrow \hat{T}$ is a finite morphism (Theorem 0.7 (i) on p. 6 in [Emer]),
\item $\lambda: \mathcal{H}(K^p)^{\text{sph}} \rightarrow \mathcal{O}(\X)$ is an $E$-algebra homomorphism,
\item $X_{cl}\subset \X(\bar{\Q}_p)$ is a Zariski-dense subset,
\end{itemize}

\noindent satisfying various properties (listed in Theorem 1.6 in [Chen], for example), the most important of which is the following: The canonical evaluation map,
$$
\text{$\X(\bar{\Q}_p)\longrightarrow (\hat{T}\times \text{Spec}  \mathcal{H}(K^p)^{\text{sph}})(\bar{\Q}_p)$, $\y$ $x \mapsto (\chi_x,\lambda_x)$,}
$$
induces a {\it{bijection}}, 
$$
X_{cl} \overset{\sim}{\longrightarrow} E(0,K^p)_{cl}.
$$
Moreover, there is a classicality criterion, analogous to Coleman's "non-critical slope implies classical", which we will not use directly (we will use that $X_{cl}$ is Zariski dense, though). More properties will be recalled below when needed, such as the connection with Emerton's Jacquet functor. 

\medskip

\noindent {\it{Notation}}. Following standard usage, by $\X(\bar{\Q}_p)$ we mean the union (or direct limit) of all $\X(L)=\Hom_{E}(\text{Sp}(L),\X)$, where $L$ ranges over all finite extensions of $E$.

\medskip

\noindent {\it{Remark}}. There is recent paper [Loef], in which Loeffler spells out how Chenevier's construction is related to Emerton's (in the case where $G(\R)$ is compact). In addition, he introduces so-called {\it{intermediate}} eigenvarieties, where one replaces $B$ with an arbitrary parabolic subgroup (and drops the assumption that $G$ should be quasi-split at $p$). It would be interesting to adapt our arguments to that setting, and thereby make progress towards the Breuil-Schneider conjecture when $\pi_p$ does not embed in a principal series (induced from the Borel). This ought to put the results of this paper, and that of [Sor], under the same roof. We hope to return to this question in the near future. 

\subsection{The Galois pseudo-character}

At each point $x \in \X(\bar{\Q}_p)$ we will assign a continuous semisimple Galois representation $\rho_x:\Gamma_{\K}\rightarrow \GL_n(\bar{\Q}_p)$, which is unramified outside 
$\Sigma=\Sigma(K^p)$, the places of $\K$ above $S(K^p)$. This is first done at a dense set of classical points, then by a formal argument one interpolates $\tr(\rho_x)$ by a pseudo-character. We refer to Chapter 1 of [BeCh] for an extensive elegant introduction to pseudo-representations, a notion going back to Wiles for $\GL(2)$, and Taylor for $\GL(n)$.

\begin{df}
Let $X_{reg}\subset X_{cl}$ be the subset of points $x$  such that $\chi_x=\psi_x\theta_x$, where $\psi_x=\otimes_{\sigma \in \Hom(F,\bar{\Q}_p)} \psi_{x,\tilde{\sigma}}$ is a regular character of $T$. That is, some $\psi_{x,\tilde{\sigma}}$ is a regular dominant character of $T_{\GL(n)}$ in the usual sense.
\end{df}

\noindent This is a Zariski-dense subset of $\X(\bar{\Q}_p)$, see p. 18 in [Chen], and the references given there. Now let $x \in X_{reg}$, and look at the corresponding pair $(\chi_x=\psi_x\theta_x,\lambda_x)$. There exists an irreducible $G(\A_f)$-summand $\pi_f$ in $\bar{\Q}_p\otimes_E H^0(\mathcal{V}_{\breve{W}_x})$, where $W_x$ has regular highest weight $\psi_x$, such that $\mathcal{H}(K^p)^{\text{sph}}$ acts on $\pi_f^{K^p}\neq 0$ via $\lambda_x$, and $\pi_p \hookrightarrow \text{Ind}_{B(\Q_p)}^{G(\Q_p)}(\theta_x)$. Thus 
$\iota^{-1}\pi_f$ is the finite part of an automorphic representation of $U(\A_F)$ of regular weight $W_x$, unramified outside $S(K^p)$, to which we can associate a continuous semisimple Galois representation
$$
\rho_x: \Gamma_{\K}\rightarrow \GL_n(\bar{\Q}_p)
$$
with the following properties:

\begin{itemize}
\item[(a)] $\rho_x^{\vee}\simeq \rho_x^c \otimes \epsilon_{cyc}^{n-1}$.
\item[(b)] For every finite place $v\nmid p$ of $F$, {\bf{outside}} $S(K^p)$, and every $w|v$ of $\K$,
the local representation $\rho_x|_{\Gamma_{\K_w}}$ is unramified, and satisfies the identity:
$$
\tr\rho_x(\text{Frob}_w)=\lambda_x(b_{w|v}(h_w)). 
$$
(Here $\text{Frob}_w$ is a geometric Frobenius, 
$h_w$ is the element of the spherical Hecke algebra for $\GL_n(\K_w)$ acting on an unramified $\Pi_w$ by $\sum \alpha_i$, where the $\alpha_i$ are the integral Satake parameters. Finally, the map
$$
b_{w|v}: \mathcal{H}(\GL_n(\K_w),K_w)\rightarrow \mathcal{H}(U(F_v),K_v)
$$
is the base change homomorphism between the spherical Hecke algebras. See [Min] for a careful useful discussion of this latter map.)
\item[(c)] For every finite place $v|p$ of $F$, the local representation $\rho_x|_{\Gamma_{\K_{\tilde{v}}}}$ is potentially semistable. Furthermore,

\begin{itemize}
\item[(i)] The semisimplification of the attached Weil-Deligne representation is
$$
WD(\rho_x|_{\Gamma_{\K_{\tilde{v}}}})^{ss}\simeq \oplus_{i=1}^n (\theta_{x,\tilde{v}}^{(i)} \circ  \text{Art}_{\K_{\tilde{v}}}^{-1}).
$$
(Here $\theta_x=\otimes_{v|p}\theta_{x,\tilde{v}}$, where $\theta_{x,\tilde{v}}$ is a smooth character of the diagonal torus $T_{\GL(n)}(\K_{\tilde{v}})\simeq (\K_{\tilde{v}}^*)^n$, factored as a product $\theta_{x,\tilde{v}}^{(1)}\otimes \cdots \otimes \theta_{x,\tilde{v}}^{(n)}$.)
\item[(ii)] The Hodge-Tate numbers are, for any embedding $\tau: \K_{\tilde{v}}\hookrightarrow \bar{\Q}_p$,
$$
\text{HT}_{\tau}(\rho_x|_{\Gamma_{\K_{\tilde{v}}}})=\{a_{\tau,j}+(n-j): j=1,\ldots, n\},
$$
where the tuple $(a_{\tau,j})$ corresponds to the dominant character $\psi_{x,v,\tau}$ of $T_{\GL(n)}$. (Here we factor $\psi_x=\otimes_{v|p}\otimes_{\tau: \K_{\tilde{v}}\hookrightarrow \bar{\Q}_p} \psi_{x,v,\tau}$.)
\end{itemize}

\end{itemize}

\noindent Observe that there may be many automorphic representations associated to a given point $x\in X_{cl}$, but they are all isomorphic outside $S(K^p)$ (and of the same weight). In particular, by (b) and Chebotarev, the Galois representation is independent of the choice of $\pi_f$, justifying the notation $\rho_x$. 

\begin{prop}
There exists a unique continuous $n$-dimensional pseudo-character $\mathcal{T}: \Gamma_{\K,\Sigma}\rightarrow \mathcal{O}(\X)^{\leq 1}$ such that
$\mathcal{T}(\text{Frob}_w)=\lambda(b_{w|v}(h_w))$ for all places $w \notin \Sigma$.
\end{prop}

\noindent {\it{Proof}}. We are in the situation of Proposition 7.1.1 in [Che]: $\X$ is reduced, $\mathcal{O}(\X)^{\leq 1}$ is a compact subring, and for all $x \in X_{reg}$, a Zariski-dense subset, we have a representation $\rho_x$ of $\Gamma_{\K,\Sigma}$ such that $\tr\rho_x(\text{Frob}_w)=\lambda(b_{w|v}(h_w))(x)$. $\square$

\begin{cor}
For every $x \in \X(\bar{\Q}_p)$, there is a unique continuous semisimple Galois representation $\rho_x:\Gamma_{\K,\Sigma}\rightarrow \GL_n(\bar{\Q}_p)$ such that
$\tr\rho_x(\text{Frob}_w)=\lambda_x(b_{w|v}(h_w))$ for all $w \notin \Sigma$.
\end{cor}

\noindent {\it{Proof}}. This follows from Theorem 1 of [Tay]. $\square$

\medskip

\noindent In particular, this applies to the classical point $x \in X_{cl}$, {\it{not}} in $X_{reg}$. One of the goals of [Chen] was to extend properties (a)-(c) above to this setting.
This was partially accomplished. See Theorem 3.3 and 3.5 of [Chen].

\section{Banach space representations}

With each point $x \in \X(L)$, we have associated an $n$-dimensional continuous pseudo-character $\mathcal{T}_x: \Gamma_{\K}\rightarrow L$, unramified outside $\Sigma(K^p)$. Here we will associate a Banach $\hat{\mathcal{H}}_L(K^p)$-module $\mathcal{B}_x$, with an admissible unitary $G(\Q_p)$-action, such that the pairs
$(\mathcal{T}_x,\mathcal{B}_x)$ form the graph of a one-to-one correspondence. We explicitly compute the locally (regular) algebraic vectors in $\mathcal{B}_x$, for $x \in X_{reg}$ such that $\mathcal{T}_x$ is absolutely irreducible, in terms of the Breuil-Schneider representation attached to $\mathcal{T}_x$, or rather, its corresponding Galois representation 
$\rho_x$. As a result, we prove the Breuil-Schneider conjecture for such $\rho_x$. 

\subsection{A global $p$-adic Langlands correspondence}

With the eigenvariety language set up, we can reformulate our findings at the end of Chapter 2. We let $X_{irr}\subset \X(\bar{\Q}_p)$ be the points $x$ for which $\rho_x$ is irreducible. 

\begin{thm}
Let $x \in X_{reg}\cap X_{irr}$, corresponding to $(\psi_x\theta_x,\lambda_x)$. Let $W_x$ be the irreducible algebraic representation of $G_{E}$ of highest weight $\psi_x$. 
Then,
$$
\bar{\Q}_p\otimes_E \tilde{H}^0(K^p)_{\text{$W_x$-alg}}^{\frak{h}=\lambda_x}\simeq (\oplus_{v|p}\widetilde{BS}(\rho_x|_{\Gamma_{\K_{\tilde{v}}}}))\otimes_{\bar{\Q}_p}
(\oplus_{\pi:\pi_{\infty}\simeq W_x,\lambda_{\pi}=\lambda_x} m_G(\pi)(\pi_f^p)^{K^p}),
$$
where we write $\frak{h}=\mathcal{H}(K^p)^{\text{sph}}$ for simplicity.
\end{thm}

\noindent This formula suggests the following definition.

\begin{df}
At each point $x \in \X(\bar{\Q}_p)$, we introduce the eigenspace
$$
\mathcal{B}_x:=(\bar{\Q}_p\otimes_E\tilde{H}^0(K^p))^{\frak{h}=\lambda_x}.
$$
This is a Banach $\hat{\mathcal{H}}(K^p)$-module with a (commuting) unitary $G(\Q_p)$-action.
\end{df}

\noindent We remind ourselves that $\mathcal{B}_x$ is nothing but the space of {\it{continuous}} $\lambda_x$-eigenforms $f:Y(K^p)\rightarrow \bar{\Q}_p$.
This sets up a one-to-one correspondence $\rho_x \leftrightarrow \mathcal{B}_x$. That is,
$$
\rho_x=\rho_{x'} \Leftrightarrow \lambda_x=\lambda_{x'} \Leftrightarrow \mathcal{B}_x=\mathcal{B}_{x'}
$$
for any two $x,x' \in \X(\bar{\Q}_p)$. Let us say that a Galois representation $\rho$ {\it{comes from}} $\X$ if 
$\rho\simeq \rho_x$ for some $x \in \X(\bar{\Q}_p)$. Similarly for Banach modules $\mathcal{B}\simeq \mathcal{B}_x$.

\medskip

\noindent This leads to the main result of this section, which in some sense is the genesis of what follows.

\begin{thm}
The eigenvariety $\X$ mediates a one-to-one correspondence between:

\begin{itemize}
\item The set of continuous semisimple Galois representations $\rho: \Gamma_{\K}\rightarrow \GL_n(\bar{\Q}_p)$ coming from $\X$. (In particular, $\rho$ is unramified outside
$\Sigma(K^p)$.)
\item The set of Banach $\hat{\mathcal{H}}(K^p)$-modules $\mathcal{B}$, with unitary $G(\Q_p)$-action, from $\X$. 
\end{itemize}

\noindent We write $\rho \leftrightarrow \mathcal{B}$ when there is a point $x \in \X(\bar{\Q}_p)$ such that $\rho\simeq \rho_x$ and $\mathcal{B}\simeq \mathcal{B}_x$.

\begin{itemize}
\item[(1)] Let $x \in \X(\bar{\Q}_p)$. If there is a regular $W$ for which $\mathcal{B}_x^{\text{$W$-alg}}\neq 0$, then $\rho_x$ is potentially semistable at all places $w|p$ of $\K$.

\item[(2)] Let $x \in X_{cl}$. Then $\mathcal{B}_x^{\text{$W_x$-alg}}\neq 0$, and $\mathcal{B}_x^{\text{$W$-alg}}=0$ for all regular $W\neq W_x$.

\item[(3)] For $x \in X_{reg}\cap X_{irr}$, the locally regular algebraic vectors of $\mathcal{B}_x$ are
$$
\mathcal{B}_x^{ralg}=\mathcal{B}_x^{\text{$W_x$-alg}}\simeq (\otimes_{v|p}\widetilde{BS}(\rho_x|_{\Gamma_{\K_{\tilde{v}}}}))\otimes_{\bar{\Q}_p}
(\oplus_{\pi:\pi_{\infty}\simeq W_x,\lambda_{\pi}=\lambda_x} m_G(\pi)(\pi_f^p)^{K^p}).
$$

\end{itemize}

\end{thm}

\noindent {\it{Proof}}. First, (1) follows from Proposition 2, which shows there is an automorphic $\pi$, with $\pi_{\infty}\simeq W$, such that $\frak{h}$ acts on $\pi_f^{K^p}$ by $\lambda_x$. Since $W$ is regular, we know how to associate a Galois representation $\rho_{\pi,\iota}$, with the usual local properties, which must be $\rho_x$ by Tchebotarev.

\medskip

\noindent For (2), we follow the same line of argument. Since $x \in X_{cl}$, there is an automorphic $\pi$ contributing to $\mathcal{B}_x^{\text{$W_x$-alg}}$. Moreover, if 
$\mathcal{B}_x^{\text{$W$-alg}}\neq 0$, there is an automorphic $\pi$, of regular weight $W$, for which $\rho_{\pi,\iota}\simeq \rho_x$. From $\rho_{\pi,\iota}$ we can recover $W$ through its Hodge-Tate numbers. Similarly for $\rho_x$, even if $x$ is not in $X_{reg}$ (this is shown in section 3.15 of [Chen], based on results of Sen, and Berger-Colmez).  Therefore, $W=W_x$. $\square$



\medskip

\noindent {\it{Remark}}. As remarked earlier, we are optimistic that one can remove the regularity hypotheses in the Theorem. Indeed it seems possible to attach Galois representations to automorphic $\pi$ of $U(\A_F)$ of {\it{irregular}} weight. When $\pi_p$ is of finite slope (that is, embeds in a principal series), this can be done by means of eigenvarieties, as in [Chen]. In general, it seems likely that one can push the ideas from the proof of Theorem 1. By [Whit], there is always a base change $\boxplus_{i=1}^t \Pi_i$, where the
$\Pi_i$ are {\it{discrete}} automorphic representations of $\GL_{n_i}(\A_{\K})$, which in turn (by the Moeglin-Waldspurger classification) are isobaric sums of cohomological, essentially conjugate self-dual, cusp forms; with which one can associate Galois representations. Local-global compatibility at $p$ follows from Caraiani's Harvard thesis. Perhaps it would be stylistically more elegant to simply admit the association of Galois representations in what follows. Then everywhere one could replace regular-algebraic vectors with algebraic vectors.

\medskip

\noindent The following result of Emerton shows how to detect the weights in $\mathcal{B}_x$. It utilizes his extension of the Jacquet functor $J_B$ to the locally analytic setting. We apply it to the locally $\Q_p$-analytic vectors $\mathcal{B}_x^{an}$.

\begin{thm}
$J_B(\mathcal{B}_x^{an})^{T(\Q_p)=\chi_x}\neq 0$, for all $x \in \X(\bar{\Q}_p)$.
\end{thm}

\noindent {\it{Proof}}. Combine Propositions 2.3.3 (iii) and 2.3.8 in [Emer]. $\square$

\subsection{Analytic variation of the Breuil-Schneider recipe}

By the very definition of $J_B$ (as an adjoint functor, see Theorem 0.3 of [Em1]), this is equivalent to having a nonzero continuous $B$-equivariant map
$$
\mathcal{C}_c^{sm}(N,\delta_B^{-1}\chi_x)\rightarrow \mathcal{B}_x^{an},
$$
where $B=TN$. If there is a so-called {\it{balanced}} map $\chi_x \hookrightarrow J_B(\mathcal{B}_x^{an})$ (Definition 0.8 in [Em2]), it is expected to arise from a nonzero continuous $G$-map (by applying $J_B$ and composing with the adjunction map)
$$
I_{\bar{B}}^G(\delta_B^{-1}\chi_x)\rightarrow \mathcal{B}_x^{an},
$$
where, loosely speaking, $I_{\bar{B}}^G$ is the closed subrepresentation of $\text{Ind}_{\bar{B}}^G$, which can be detected by its Jacquet module (Corollary 5.1.4 in [Em2]). This expectation is known in many cases, see Theorem 0.13 in [Em2]. Note that, for $x \in X_{irr}$,
$$
I_{\bar{B}}^G(\delta_B^{-1}\chi_x)\twoheadrightarrow \widetilde{BS}(\rho_x):=\otimes_{v|p}\widetilde{BS}(\rho_x|_{\Gamma_{\K_{\tilde{v}}}})
$$
according to Remark 5.1.8 in [Em2] (an isomorphism when $\text{Ind}_B^G(\theta_x)$ is irreducible). What we wish to point out in this paragraph, is that $I_{\bar{B}}^G(\delta_B^{-1}\chi_x)$ {\it{always}} satisfies the Emerton condition (that is, (3) in the introduction), for every point $x$ on the eigenvariety $\X$, even though we do {\it{not}} know if it admits a $G$-invariant norm: 

\begin{prop}
$I_{\bar{B}}^G(\delta_B^{-1}\chi_x)$ satisfies the Emerton condition for all $x \in \X(\bar{\Q}_p)$.
\end{prop}

\noindent {\it{Proof}}. Since $\mathcal{B}_x^{an}$ {\it{has}} an invariant norm, it satisfies the Emerton condition. (See Lemma 4.4.2 in [Em1].) In particular, 
$|\delta_B^{-1}(z)\chi_x(z)|_p\leq 1$ for $z \in T^+$, since $\chi_x \hookrightarrow J_B(\mathcal{B}_x^{an})$. Therefore, the sup-norm $\|\cdot\|$on $\mathcal{C}_c^{sm}(N,\delta_B^{-1}\chi_x)$ is non-increasing under the action of the submonoid $NT^+$. By restriction to $N$,
$$
I_{\bar{B}}^G(\delta_B^{-1}\chi_x)^{\frak{n}}:=\cup_{N_0}I_{\bar{B}}^G(\delta_B^{-1}\chi_x)^{N_0}\hookrightarrow \mathcal{C}_c^{sm}(N,\delta_B^{-1}\chi_x).
$$
Consequently $I_{\bar{B}}^G(\delta_B^{-1}\chi_x)^{\frak{n}}$ carries a norm such that $\|gx\|\leq \|x\|$ for all $g \in NT^+$. This is sufficient to run the argument proving Lemma 4.4.2 in [Em1]. $\square$

\section{Weak local-global compatibility}

In this section we deduce from our previous results that $\widetilde{BS}(\rho_x)$ admits an invariant norm such that the completion satisfies (a strong version of) local-global compatibility. However, we cannot show that this completion  $\widetilde{BS}(\rho_x)^{\wedge}$ only depends on the restrictions of $\rho_x$ to places above $p$.
Ultimately, we will restrict ourselves to the {\it{unramified}} case, and prove a weak version of local-global compatibility (somewhat similar to part (1) of Theorem 1.2.1 in [Eme]). The $p$-adic local Langlands correspondence, still mysterious in higher rank, is replaced by the coarse version in [ScTe], which associates a huge Banach representation
$B_{\xi,\zeta}$ with a pair $(\xi,\zeta)$ satisfying the Emerton condition (here $\xi$ is an irreducible algebraic representation, and $\zeta$ is a suitable Weyl-orbit in the dual torus).
The philosophy propounded in [ScTe] and [BrSc] is that the completed quotients of  $B_{\xi,\zeta}$ should somehow correspond to the crystalline representations of type  $(\xi,\zeta)$. This is well-understood for $\GL_2(\Q_p)$, where the admissible filtration is usually unique, and $B_{\xi,\zeta}$ essentially {\it{is}} the local $p$-adic Langlands correspondence in the crystalline case. We provide evidence supporting this philosophy for $n>2$.

\subsection{Completions of the algebraic vectors}

{\it{Split ramification and the automorphic representation $\pi_x$}}. Throughout, we will make the assumption that we have {\it{split ramification}}. That is,
$S(K^p)\subset \text{Spl}_{\K|F}$. This has the effect that local base change $\text{BC}_{w|v}$ is defined at {\it{all}} places $v$. We fix a point $x \in X_{reg}\cap X_{irr}$, as above.
Under our ramification hypothesis, there is a {\it{unique}} automorphic $\pi$ contributing to the (regular) algebraic vectors $\mathcal{B}_x^{ralg}$ in Theorem 4, part (3). Indeed,
any such $\pi$ has an irreducible Galois representation $\rho_{\pi,\iota}\simeq \rho_x$, and therefore $\text{BC}_{\K|F}(\pi)$ must be cuspidal, and it is uniquely determined at the infinite places, and away from $\Sigma(K^p)$. By strong multiplicity one for $\GL_n$, the base change is unique. Locally, $\text{BC}_{w|v}$ is injective (see Corollary 4.2 in [Min]), and therefore $\pi$ is uniquely determined. We denote it $\pi_x=\otimes \pi_{x,v}$. Its local components $\pi_{x,v}$ are given by
$$
WD(\rho_x|_{\Gamma_{\K_w}})^{F-ss}\simeq rec(BC_{w|v}(\pi_{x,v})\otimes |\det|_w^{(1-n)/2}).
$$
We think of $\{\pi_x\}$ as a family of automorphic representations interpolated by $\X$. In general (without split ramification) the $\pi_x$ will be $L$-packets, not singletons.

\medskip

\noindent With this notation, part (3) of Theorem 4 becomes: For all $x \in X_{reg}\cap X_{irr}$,
$$
\mathcal{B}_x^{ralg}\simeq \widetilde{BS}(\rho_x) \otimes (\otimes_{v\nmid p}\pi_{x,v}^{K_v})^{m(\pi_x)}.
$$
Most likely, $m(\pi_x)=1$, and this may already be in the literature. However, we have not been able to find a suitable reference. Now, since 
$\otimes_{v\nmid p}\pi_{x,v}^{K_v}$ is a simple $\mathcal{H}(K^p)$-module, we may think of $\widetilde{BS}(\rho_x)^{m(\pi_x)}$ as its multiplicity space in 
$\mathcal{B}_x^{ralg}$,
$$
\widetilde{BS}(\rho_x)^{m(\pi_x)} \overset{\sim}{\longrightarrow} \Hom_{\mathcal{H}(K^p)}(\otimes_{v\nmid p}\pi_{x,v}^{K_v}, \mathcal{B}_x^{ralg}),
$$
as representations of $G(\Q_p)$. We will view the right-hand-side as sitting inside a Banach space of continuous transformations. For that purpose, we first look at each local component $\pi_{x,v}$, where $v \nmid p$. When $v$ splits, it can be identified with a $p$-integral irreducible representation of $\GL_n(F_v)$. By Theorem 1 in [Vig], it has a
unique commensurability class of stable lattices. Correspondingly, $\pi_{x,v}$ has a unique equivalence class of $\GL_n(F_v)$-invariant norms $\|\cdot\|_v$. (By Theorem 1 in [Vign], the completion $\hat{\pi}_{x,v}$ is a topologically irreducible unitary Banach space representation of $\GL_n(F_v)$.) When $\pi_{x,v}$ is unramified, its Satake parameters are $p$-units, and one easily finds a stable lattice in a suitable unramified principal series, again resulting in a $U(F_v)$-invariant (supremum) norm $\|\cdot\|_v$, which we may normalize such that a given spherical vector has norm one. The tensor product norm (see Proposition 17.4 in [Sc]) on $\otimes_{v\nmid p}\pi_{x,v}$ is then invariant under $G(\A_f^p)$.
 By restriction, the finite-dimensional space $\otimes_{v\nmid p}\pi_{x,v}^{K_v}$ inherits a norm, and becomes a Banach-module for $\hat{\mathcal{H}}(K^p)$. With this extra structure at hand, 
$$
\Hom_{\mathcal{H}(K^p)}(\otimes_{v\nmid p}\pi_{x,v}^{K_v}, \mathcal{B}_x^{ralg})\hookrightarrow \mathcal{L}_{\hat{\mathcal{H}}(K^p)}(\otimes_{v\nmid p}\pi_{x,v}^{K_v},\mathcal{B}_x).
$$
(Here $\mathcal{L}$ denotes the space of continuous linear transformations, equipped with the usual transformation norm, see Corollary 3.2 in [Sc].) We have to check that any $\mathcal{H}(K^p)$-equivariant map $\otimes_{v\nmid p}\pi_{x,v}^{K_v} \overset{\phi}{\rightarrow} \mathcal{B}_x$ is automatically continuous: If $\phi\neq 0$, it must be injective (by simplicity), and thus $\|\phi(\cdot)\|_{\mathcal{B}_x}$ defines a norm on $\otimes_{v\nmid p}\pi_{x,v}^{K_v}$. However, all norms on a finite-dimensional space are equivalent (4.13 in [Sc]), so
that $\|\phi(u)\|_{\mathcal{B}_x}\leq C\|u\|$ for some constant $C>0$, and all $u$. Altogether, this embeds $\widetilde{BS}(\rho_x)$ into a Banach space (Proposition 3.3 in [Sc]):
\begin{equation}
\widetilde{BS}(\rho_x)^{m(\pi_x)} \hookrightarrow \mathcal{L}_{\hat{\mathcal{H}}(K^p)}(\otimes_{v\nmid p}\pi_{x,v}^{K_v},\mathcal{B}_x).
\end{equation}
If we restrict the tranformation norm to $\widetilde{BS}(\rho_x)^{m(\pi_x)}$, we arrive at:

\begin{cor}
Let $x \in X_{reg}\cap X_{irr}$ be a point such that $m(\pi_x)=1$. Then there is a $G(\Q_p)$-invariant norm $\|\cdot\|$ on $\widetilde{BS}(\rho_x)$ such that the corresponding completion $\widetilde{BS}(\rho_x)^{\wedge}$ satisfies the following: There is a topological isomorphism,
$$
\widetilde{BS}(\rho_x)^{\wedge}\otimes (\otimes_{v\nmid p}\pi_{x,v}^{K_v}) \overset{\sim}{\longrightarrow} \overline{\mathcal{B}_x^{ralg}},
$$
where $\overline{\mathcal{B}_x^{ralg}}$ is the closure of the regular-algebraic vectors $\mathcal{B}_x^{ralg}$ in $\mathcal{B}_x$. Moreover,
\begin{itemize}
\item $\widetilde{BS}(\rho_x)^{\wedge}$ is an \underline{admissible} unitary Banach space representation of $G(\Q_p)$.
\item Its regular-algebraic vectors $\widetilde{BS}(\rho_x)$ form a \underline{dense} subspace.
\end{itemize}
\end{cor}

\noindent {\it{Proof}}. We obtain $\|\cdot\|$ by restricting the transformation norm to $\widetilde{BS}(\rho_x)$. Thus (2) becomes an isometry, and extends uniquely to an isometry
$$
\widetilde{BS}(\rho_x)^{\wedge} \hookrightarrow \mathcal{L}_{\hat{\mathcal{H}}(K^p)}(\otimes_{v\nmid p}\pi_{x,v}^{K_v},\mathcal{B}_x).
$$
To ease the notation, let us write $M=\otimes_{v\nmid p}\pi_{x,v}^{K_v}$ throughout this proof; a finite-dimensional simple $\mathcal{H}(K^p)$-module. We tensor the isometry by this $M$,
$$
j:\widetilde{BS}(\rho_x)^{\wedge}\otimes M \hookrightarrow \mathcal{L}_{\hat{\mathcal{H}}(K^p)}(M,\mathcal{B}_x)\otimes M\overset{\sim}{\longrightarrow} \mathcal{B}_x[M].
$$
Here $\mathcal{B}_x[M]$ denotes the closure of the sum of all closed $\mathcal{H}(K^p)$-submodules of $\mathcal{B}_x$ isomorphic to $M$ (a topological direct sum of a subcollection, by Zorn). Note that $\End_{\mathcal{H}(K^p)}(M)=\bar{\Q}_p$. Note also that the tensor products (equipped with their tensor product norms, as on p. 110 in [Sc]) are already complete, as $M$ is finite-dimensional. The above isomorphism with $\mathcal{B}_x[M]$ is a {\it{topological}} isomorphism by the open mapping theorem (8.6, p. 55 in [Sc]), but not necessarily isometric. Consequently, $\im(j)\subset \mathcal{B}_x$ is a closed subspace, containing $\mathcal{B}_x^{ralg}$ by Theorem 4. In fact, $\im(j)$ is the closure of $\mathcal{B}_x^{ralg}$ in $\mathcal{B}_x$, since $\widetilde{BS}(\rho_x)$ is dense in the completion $\widetilde{BS}(\rho_x)^{\wedge}$. Again invoke the open mapping theorem to see that $j$ is a topological isomorphism onto $\overline{\mathcal{B}_x^{ralg}}$. Admissibility of $\widetilde{BS}(\rho_x)^{\wedge}$ follows from admissibility of $\mathcal{B}_x$. $\square$

\medskip

\noindent {\it{Remark}}.  Equivalently, there is a $G(\Q_p)$-equivariant topological isomorphism,
$$
\widetilde{BS}(\rho_x)^{\wedge} \overset{\sim}{\longrightarrow} \mathcal{L}_{\hat{\mathcal{H}}(K^p)}(\otimes_{v\nmid p}\pi_{x,v}^{K_v},\overline{\mathcal{B}_x^{ralg}}).
$$
We like to think of this Banach space representation $\widetilde{BS}(\rho_x)^{\wedge}$ as a "rough" candidate for a $p$-adic local Langlands correspondence
$\frak{B}(\rho_x)=\hat{\otimes}_{v|p} \frak{B}(\rho_x|_{\Gamma_{\K_{\tilde{v}}}})$, at least when the various restrictions $\rho_x|_{\Gamma_{\K_{\tilde{v}}}}$ are irreducible. Of course, to really justify this point of view, one would need to show that 
the completion $\widetilde{BS}(\rho_x)^{\wedge}$ only depends on the restrictions $\rho_x|_{\Gamma_{\K_{\tilde{v}}}}$ at $p$, and that it factors as a tensor product $\hat{\otimes}_{v|p}$ of appropriate completions $\widetilde{BS}(\rho_x|_{\Gamma_{\K_{\tilde{v}}}})^{\wedge}$. Both appear to be very difficult questions.

\subsection{Universal modules: The crystalline case}

We now specialize to the {\it{crystalline}} case, where we can relate $\widetilde{BS}(\rho_x)^{\wedge}$ to the Schneider-Teitelbaum universal modules $B_{\xi,\zeta}$, which are given by a purely local construction at $p$. They are expected to be quite large. However, for $n>2$ it is not even known that $B_{\xi,\zeta}\neq 0$ (Conjecture 6.1, p. 24 in [BrSc]).
For $n=2$ this is a deep result of Berger and Breuil. We will prove non-vanishing when $(\xi,\zeta)$ "comes from an eigenvariety". This will be a by-product of a stronger result.

\begin{df}
A classical point $x \in X_{cl}$ is called \underline{old} if $\rho_x$ is crystalline at all places above $p$. That is, $\Hom_{G(\Z_p)}(W_x,\mathcal{B}_x)\neq 0$.
Equivalently, $\pi_{x,v}$ is unramified for all $v|p$. We denote the set of old points by $X_{old}$.
\end{df}

\noindent Thus, from now on, we fix a point $x \in X_{reg}\cap X_{irr} \cap X_{old}$. By Proposition 1,
$$
 \widetilde{BS}(\rho_x)=W_x \otimes \pi_{x,p}\overset{\sim}{\longrightarrow} W_x \otimes \text{Ind}_{B}^{G}(\theta_x),
$$
where $\theta_x$ is unramified smooth. (Indeed, for any point $x$, $\pi_{x,p}$ embeds into the (non-normalized) principal series $\text{Ind}_B^G(\theta_x)$. Since $x$ is old, 
$\pi_{x,p}$ is unramified, and hence so is $\theta_x$. Furthermore, as $x \in X_{irr}$, the base change $\text{BC}_{\K|F}(\pi_x)$ is cuspidal, and therefore generic. In particular, 
$\pi_{x,p}$ must be generic. As is well-known, this implies that $\pi_{x,p}$ must be the full unramified principal series.)

\medskip

\noindent As in [ScTe] we express $\pi_{x,p}\simeq \text{Ind}_B^G(\theta_x)$ in terms of the {\it{universal module}}. This goes back to Borel and Matsumoto, and is defined as follows. For any algebra character $\zeta: \mathcal{H}(G,K) \rightarrow \bar{\Q}_p$
(where $K=G(\Z_p)$ is hyperspecial when $p$ is assumed to be unramified in $F$) we introduce the smooth representation
$$
\mathcal{M}_{\zeta}=\text{c-Ind}_{K}^{G}(1)\otimes_{\mathcal{H}(G,K),\zeta} \bar{\Q}_p=\mathcal{C}_c(K\backslash G,\bar{\Q}_p)\otimes_{\mathcal{H}(G,K),\zeta} \bar{\Q}_p.
$$ 
The pair $(\mathcal{M}_{\zeta},1_K)$ is a universal initial object in the category of pairs $(V,v)$, where $V$ is an unramified smooth representation of $G(\Q_p)$, and $v \in V^K$ is a nonzero vector on which $\mathcal{H}(G,K)$ acts via $\zeta$. That is, there is a unique $G(\Q_p)$-map
$\mathcal{M}_{\zeta} \rightarrow V$ which maps $1_K\mapsto v$. The image of this map is the span of the orbit $Gv$ (since $\mathcal{M}_{\zeta}$ is generated by $1_K$). In what follows we will take $\zeta_x=\hat{\theta}_x$, the eigensystem of $\text{Ind}_B^G(\theta_x)^K$. The choice of a spherical vector yields
$$
\text{$\mathcal{M}_{\zeta_x}\rightarrow \text{Ind}_B^G(\theta_x)$, $\y$ $\zeta_x=\hat{\theta}_x$.}
$$
It is a general fact that the two representations have the same semi-simplification (see the Orsay Ph.D. thesis of X. Lazarus for a thorough discussion in greater generality). Under our assumptions, $\text{Ind}_B^G(\theta_x)$ is irreducible, and therefore the above must be an isomorphism. Consequently, we may identify 
$$
\widetilde{BS}(\rho_x)\simeq W_x \otimes \mathcal{M}_{\zeta_x}\simeq \text{c-Ind}_{K}^{G}(\xi_x)\otimes_{\mathcal{H}_{\xi_x}(G,K),\zeta_x} \bar{\Q}_p=:H_{\xi_x,\zeta_x}.
$$
Here we have changed notation $\xi_x:=W_x$ to aid comparison with [ScTe]. The algebra $\mathcal{H}_{\xi_x}(G,K)$ is by definition the $G$-endomorphisms of 
$\text{c-Ind}_{K}^{G}(\xi_x)$. Or, more concretely, compactly supported $K$-biequivariant functions $G \rightarrow \End(\xi_x)$, with convolution. However, since $\xi$ is an irreducible representation of $G$ (viewed as a representation of $K$), as on p. 639 in [ScTe] one can identify the algebras
$$
\text{$\mathcal{H}(G,K)\overset{\sim}{\longrightarrow} \mathcal{H}_{\xi_x}(G,K)$, $\y$ $h \mapsto (g \mapsto h(g)\xi_x(g))$.}
$$
In the definition of $H_{\xi_x,\zeta_x}$ we view $\zeta_x$ as a character of $\mathcal{H}_{\xi_x}(G,K)$ via this isomorphism, as 
at the bottom of p. 670 in [ScTe], where $H_{\xi,\zeta}$ is defined.

\medskip

\noindent The representation $H_{\xi_x,\zeta_x}$ has a natural locally convex topology, being a quotient of $\text{c-Ind}_{K}^{G}(\xi_x)$, which has a supremum norm: Pick any norm $\|\cdot\|_{\xi_x}$ on $\xi_x$, which is invariant under (the compact group) $K$. They are all equivalent since $\xi_x$ is finite-dimensional (4.13 in [Sc]). As in Appendix A, for $f \in \text{c-Ind}_{K}^{G}(\xi_x)$, we let
$$
\|f\|_{\xi_x,\infty}={\sup}_{g \in G(\Q_p)} \|f(g)\|_{\xi_x}<\infty
$$
This defines an norm $\|\cdot\|_{\xi_x,\infty}$ on the compact induction, which is obviously invariant under $G(\Q_p)$, and it induces a quotient {\it{seminorm}} on the representation
$$
H_{\xi_x,\zeta_x}=(\text{c-Ind}_{K}^{G}(\xi_x))/(\ker\zeta_x)(\text{c-Ind}_{K}^{G}(\xi_x)).
$$
We will show below that in fact this is a {\it{norm}}, but this is far from clear a priori!

\medskip

\noindent Following [ScTe], on p. 671 where they define $B_{\xi,\zeta}$, we introduce the space
$$
B_{\xi_x,\zeta_x}:=\hat{H}_{\xi_x,\zeta_x}=(H_{\xi_x,\zeta_x}/\overline{\{0\}})^{\wedge}=\text{Hausdorff completion of $H_{\xi_x,\zeta_x}$}.
$$
(We refer to 7.5 in [Sc] for a general discussion of Hausdorff completions.) We have defined a Banach space $B_{\xi_x,\zeta_x}$ with a unitary $G(\Q_p)$-action
However, it is not clear at all that it is {\it{nonzero}}. This is in fact a fundamental problem! Conjecture 6.1 on p. 24 in [BrSc] says that $B_{\xi,\zeta}\neq 0$ whenever the Emerton condition is satisfied (the converse is known). This follows from our methods when the pair $(\xi,\zeta)$ comes from an eigenvariety. That is, when it is of the form 
$(\xi_x,\zeta_x)$ for an old irreducible point $x$. What we prove is a strengthening:

\begin{thm}
Let $x \in X_{reg}\cap X_{irr}\cap X_{old}$ be a classical point such that $m(\pi_x)=1$. Then $H_{\xi_x,\zeta_x}$ is Hausdorff, $B_{\xi_x,\zeta_x} \neq 0$ is its universal completion. Furthermore:
\begin{itemize}
\item[(a)] There is a continuous map, with dense image, $B_{\xi_x,\zeta_x} \rightarrow \widetilde{BS}(\rho_x)^{\wedge}$
(into the completion from Corollary 3) which restricts to an isomorphism $H_{\xi_x,\zeta_x}\overset{\sim}{\longrightarrow} \widetilde{BS}(\rho_x)$
onto the regular-algebraic vectors.
\item[(b)] There is a \underline{nonzero} $G(\Q_p)\times \hat{\mathcal{H}}(K^p)$-equivariant continuous map
$$
B_{\xi_x,\zeta_x} \otimes (\otimes_{v\nmid p}\pi_{x,v}^{K_v}) \rightarrow \overline{\mathcal{B}_x^{ralg}}
$$ 
with dense image.
\end{itemize}
\end{thm}

\noindent {\it{Proof}}. From (2) of the previous section, we have a $G(\Q_p)$-embedding
$$
H_{\xi_x,\zeta_x}\simeq \widetilde{BS}(\rho_x) \hookrightarrow \mathcal{L}_{\hat{\mathcal{H}}(K^p)}(M,\mathcal{B}_x),
$$
where we keep writing $M=\otimes_{v\nmid p}\pi_{x,v}^{K_v}$. We claim this map is automatically continuous, when we equip the $\mathcal{L}$-space with the transformation norm, and $H_{\xi_x,\zeta_x}$ with the quotient seminorm induced by $\|\cdot\|_{\xi_x,\infty}$. Since $H_{\xi_x,\zeta_x}$ gets the quotient topology, we just have to check continuity of the inflated map
$$
\text{c-Ind}_{K}^{G}(\xi_x) \twoheadrightarrow H_{\xi_x,\zeta_x} \hookrightarrow \mathcal{L}_{\hat{\mathcal{H}}(K^p)}(M,\mathcal{B}_x).
$$
This is part of Lemma 3 in Appendix A (in which Frobenius reciprocity is made explicit). In particular, the seminorm on $H_{\xi_x,\zeta_x}$ is actually a norm (as the kernel of the above map is closed). Therefore, $H_{\xi_x,\zeta_x}$ is Hausdorff, and $B_{\xi_x,\zeta_x}$ is its universal completion.  That is, there is an isometry with dense image,
$$
H_{\xi_x,\zeta_x} \hookrightarrow B_{\xi_x,\zeta_x}.
$$
($\Rightarrow B_{\xi_x,\zeta_x}$ is nonzero.) By continuity of the initial map, it has a unique extension
$$
B_{\xi_x,\zeta_x} \rightarrow \mathcal{L}_{\hat{\mathcal{H}}(K^p)}(M,\mathcal{B}_x),
$$
which is continuous (but not necessarily injective) and maps into the completion $\widetilde{BS}(\rho_x)^{\wedge}$ from Corollary 3, with dense image (but not necessarily onto). 
$\square$

\medskip

\noindent {\it{Remark}}. This fits perfectly with the picture suggested in the papers [ScTe] and [BrSc]. If there is a local $p$-adic Langlands correspondence $\rho \mapsto \frak{B}(\rho)$, these references speculate that $B_{\xi,\zeta}$ maps to each $\frak{B}(\rho)$, with dense image, for all crystalline representations $\rho$ of type $(\xi,\zeta)$. The previous Theorem provides strong evidence that $\widetilde{BS}(\rho_x)^{\wedge}$ is at least closely related to the (elusive) local $p$-adic Langlands correspondence $\frak{B}(\rho_x)=\hat{\otimes}_{v|p} \frak{B}(\rho_x|_{\Gamma_{\K_{\tilde{v}}}})$.


\section{On the density of the algebraic vectors}

In general, it is {\it{not}} expected that $\mathcal{B}_x^{alg}$ is dense in $\mathcal{B}_x$. In this section, we will adapt (and elaborate on) an argument from sections 5.3 and 5.4 in [Eme], which shows that the algebraic vectors {\it{are}} dense, if we substitute $\mathcal{B}_x$ by a bigger space.

\subsection{Injectivity of certain modules}

We fix a finite extension $L|\Q_p$, and we will write $\mathcal{O}=\mathcal{O}_L$ and $\varpi=\varpi_L$, and so on. We will look at locally constant functions $f:Y(K^p)\rightarrow A$, taking values in various finite $\mathcal{O}$-modules $A=\mathcal{O}/\varpi^s\mathcal{O}$, where $s$ is a positive integer.
These functions form a (discrete) torsion $\mathcal{O}$-module, denoted $H^0(K^p,A)$, carrying a natural action of the Hecke algebra $\mathcal{H}_{\mathcal{O}}(K^p)$, and a commuting smooth $G(\Q_p)$-action, which is {\it{admissible}} in the following sense: For every compact open subgroup of $G(\Q_p)$, its invariants form a finite $\mathcal{O}$-module (torsion and finitely generated means finite cardinality, since $A$ is a finite ring).

\begin{lem}
Suppose $K^p$ is sufficiently small (for example, it suffices that $K_v$ has no $p$-torsion for some $v\nmid p$). Then, for any compact open subgroup $K_p \subset G(\Q_p)$,
$$
\text{$H^0(K^p,\mathcal{O}/\varpi^s\mathcal{O})$ is an {\bf{injective}} smooth $(\mathcal{O}/\varpi^s\mathcal{O})[K_p]$-module for all $s \geq 1$.}
$$
Consequently, every direct summand of $H^0(K^p,\mathcal{O}/\varpi^s\mathcal{O})$ is an injective module.
\end{lem}

\noindent {\it{Proof}}. We have to show the exactness of the functor sending a module $M$ to 
$$
\Hom_{\mathcal{O}[K_p]}(M,H^0(K^p,\mathcal{O}/\varpi^s\mathcal{O})).
$$
Here $M$ is an $\mathcal{O}[K_p]$-module with $\varpi^sM=0$. Therefore, it has Pontryagin dual,
$$
M^{\vee}=\Hom_{\mathcal{O}}(M,L/\mathcal{O})=\Hom_{\mathcal{O}}(M,\varpi^{-s}\mathcal{O}/\mathcal{O})\simeq \Hom_{\mathcal{O}/\varpi^s\mathcal{O}}(M,\mathcal{O}/\varpi^s\mathcal{O}).
$$
(Here $M$ is smooth, so we equip it with the discrete topology.) The initial module above can then be identified with that consisting of all functions,
$$
\text{$f: Y(K^p) \rightarrow M^{\vee}$, $\y$ $f(gk)=k^{-1}f(g)$,}
$$
for $k \in K_p$. Choosing representatives $g_i \in G(\A_f)$ for the finite set $Y(K_pK^p)$,
and mapping $f$ to the tuple of all $f(g_i)$, then identifies the latter with the direct sum $\oplus_i (M^{\vee})^{\Gamma_i}$, where the $\Gamma_i$ are certain finite subgroups of $K_p$, having prime-to-$p$ order by assumption. This ensures that $(\cdot)^{\Gamma_i}$ is exact, by averaging.
Also, taking the Pontryagin dual is exact ($L/\mathcal{O}$ is divisible). Finally, as is well known (and easy to check) every summand of an injective module is itself injective. $\square$

\medskip

\noindent {\it{Examples}}. Let us first introduce certain finite type Hecke algebras. For each $K_p$, we let $\T(K_pK^p)$ denote the image of $\frak{h}^{\circ}=\mathcal{H}_{\mathcal{O}}(K^p)^{\text{sph}}$ in the endomorphism algebra $\End_{\mathcal{O}}H^0(Y(K_pK^p),\mathcal{O})$. Thus $\T(K_pK^p)$ is finite free over (the PID) $\mathcal{O}$, and we endow it with the $\varpi$-adic topology. If we have a subgroup $K_p'\subset K_p$, there is a natural restriction map $\T(K_p'K^p)\rightarrow \T(K_pK^p)$, and we take the limit,
$$
\T(K^p):=\underset{K_p}{\varprojlim} \T(K_pK^p)\subset \End_{\mathcal{O}}\tilde{H}^0(K^p)^{\circ},
$$
the closure of the image of $\mathcal{H}_{\mathcal{O}}(K^p)^{\text{sph}}$. This defines a reduced, commutative, complete, topological $\mathcal{O}$-algebra. Moreover, 
$\T(K^p)$ has only {\it{finitely}} many maximal ideals: They correspond to the maximal ideals of $\T(K^p)\otimes \F$, which is the image of $\frak{h}^{\circ}$ in
$\End_{\F}H^0(K^p,\F)$. Hence, the maximal ideals are in bijection with the (Galois conjugacy classes of) eigensystems $\frak{h}^{\circ} \rightarrow \F$ which occur in 
$H^0(K^p,\F)$. If $K_p$ is any pro-$p$ group, they all must occur in $H^0(Y(K_pK^p),\F)$, which is finite-dimensional. Therefore, since $\mathcal{O}$ is complete, we have
$$
\T(K^p) \overset{\sim}{\longrightarrow} \prod_{\m}\T(K^p)_{\m},
$$
where the product extends over the finitely many maximal ideals $\m \subset \T(K^p)$, and $\T(K^p)_{\m}$ denotes the corresponding localization, a complete local $\mathcal{O}$-algebra. (We refer to Chapter 4 of [DDT] for a discussion of the commutative algebra needed.) We will use this product decomposition as follows:
Obviously,
$$
\tilde{H}^0(K^p)^{\circ}/\varpi^s \tilde{H}^0(K^p)^{\circ}\simeq H^0(K^p,\mathcal{O}/\varpi^s\mathcal{O})
$$
carries an action of $\T(K^p)$. This gives rise to a {\it{direct sum}} decomposition,
$$
H^0(K^p,\mathcal{O}/\varpi^s\mathcal{O}) \overset{\sim}{\longrightarrow} \bigoplus_{\m} H^0(K^p,\mathcal{O}/\varpi^s\mathcal{O})_{\m},
$$
into localized smooth admissible $G(\Q_p)$-submodules over $\mathcal{O}/\varpi^s\mathcal{O}$,
$$
H^0(K^p,\mathcal{O}/\varpi^s\mathcal{O})_{\m}:=H^0(K^p,\mathcal{O}/\varpi^s\mathcal{O})\otimes_{\T(K^p)}\T(K^p)_{\m},
$$
which are then {\it{injective}} $(\mathcal{O}/\varpi^s\mathcal{O})[K_p]$-modules, for every compact open $K_p$. 

\medskip

\noindent To connect this to the previous discussion, one could take the maximal ideal $\m_x=\ker(\bar{\lambda}_x)$ for a point $x \in \X(L)$. A priori this is a maximal ideal in $\frak{h}^{\circ}$, but it is the pull-back of an ideal $\m\subset \T(K^p)$ since $\bar{\lambda}_x$ occurs in tame level $K^p$.

\subsection{Projective modules over certain Iwasawa algebras}

To simplify notation, we write $A=\mathcal{O}/\varpi^s\mathcal{O}$ in this section, where $s>0$ is fixed for the moment. We will briefly recall known facts about the Pontryagin duality functor $M \mapsto M^{\vee}$, which sends a discrete $A[K_p]$-module $M$ to the compact
$$
M^{\vee}=\Hom_{\mathcal{O}}(M,L/\mathcal{O})\simeq \Hom_A(M,A).
$$
If $M$ is smooth, $M=\underset{H}{\varinjlim}M^H$, with $H$ running over normal open subgroups of $K_p$, and therefore its dual $M^{\vee}=\underset{H}{\varprojlim}(M^H)^{\vee}$ becomes a module for 
$$
A[[K_p]]:=\underset{H}{\varprojlim} A[K_p/H],
$$
the {\it{Iwasawa algebra}}. Conversely, if $X$ is an $A[[K_p]]$-module, $X/I_HX$ becomes a module for $A[K_p/H]$, where $I_H$ is the kernel of the natural projection $A[[K_p]]\twoheadrightarrow A[K_p/H]$. It follows that $X^{\vee}$ is again a smooth $A[K_p]$-module, since
$$
(X^{\vee})^H \simeq (X/I_HX)^{\vee}.
$$
Thus, duality sets up a one-to-one correspondence $M \leftrightarrow X$ between smooth discrete $A[K_p]$-modules and compact $A[[K_p]]$-modules, which reverses arrows.

\begin{lem}
Suppose $M$ is a smooth $A[K_p]$-module, with Pontryagin dual $M^{\vee}$.
\begin{itemize}
\item[(i)] $M$ is admissible $\Longleftrightarrow$ $M^{\vee}$ is finitely generated over $A[[K_p]]$.
\item[(ii)] $M$ is injective $\Longleftrightarrow$ $M^{\vee}$ is a projective $A[[K_p]]$-module.
\end{itemize}
\end{lem}

\noindent {\it{Proof}}. For part (i), if $X$ is finitely generated over $A[[K_p]]$, we deduce that $X/I_HX$ is finitely generated over $A[K_p/H]$, which is a ring of finite cardinality. Therefore its dual $M^H$ is (physically) finite. For the converse, suppose $M$ is admissible. Then, first of all, $M^{\vee}$ is profinite, so we may apply the "converse" (topological) Nakayama lemma discussed in depth in [BH] (specifically, their main Theorem in Chapter 3, section (1), and its Corollary): To verify that $M^{\vee}$ is finitely generated over the compact ring $A[[K_p]]$, it suffices to check that $X/I_H X$ is finitely generated over $A[K_p/H]$, for {\it{some}} $H$ such that $I_H^n \rightarrow 0$ as $n \rightarrow \infty$.
This limit holds for any pro-$p$-group $H$, see Lemma 3.2 in [ScT], for example. Finiteness of $X/I_H X$, or rather its dual $M^H$, is admissibility.

\medskip

\noindent For part (ii), use that Pontryagin duality is exact (divisibility of $L/\mathcal{O}$). It follows that $\Hom_{A[K_p]}(-,M)$ is exact if and only if $\Hom_{A[[K_p]]}(M^{\vee},-)$ is exact. $\square$

\medskip

\noindent From the last two lemmas, we immediately conclude the following:

\begin{prop}
Suppose $K^p$ is sufficiently small. Then, for any compact open subgroup $K_p \subset G(\Q_p)$, the dual 
$H^0(K^p,A)^{\vee}$ is a {\bf{projective}} finitely generated module over $A[[K_p]]$ for all $s \geq 1$.
The same is true for any direct summand, such as the localized module $H^0(K^p,A)_{\m}^{\vee}$ for any maximal ideal $\m$.
\end{prop}

\noindent For later use, we will record the following fact here. Often, the Iwasawa algebra $A[[K_p]]$ is viewed as a distribution algebra. Indeed, there is a natural pairing with the continuous (that is, locally constant) functions $\mathcal{C}(K_p,A)$.

\begin{lem}
$A[[K_p]]\overset{\sim}{\longrightarrow}\mathcal{C}(K_p,A)^{\vee}$, as modules over $A[[K_p]]$.
\end{lem}

\noindent {\it{Proof}}. For any normal open subgroup $H$, there is a canonical integration pairing,
$$
\text{$\mathcal{C}(K_p/H,A)\times A[K_p/H]\rightarrow A$, $\y$ $(f,\mu)\mapsto {\sum}_{k \in K_p/H}f(k)\mu(k)$,}
$$
which is non-degenerate, and therefore defines an isomorphism 
$$
\text{$A[K_p/H]\overset{\sim}{\longrightarrow}\mathcal{C}(K_p/H,A)^{\vee}$, $\y$ $\mu \mapsto (-,\mu)$.}
$$
This is easily checked to preserve the $A[K_p/H]$-module structures on both sides. Moreover, as $H$ varies, these isomorphisms are compatible with the transition maps. Passing to the projective limit $\underset{H}{\varprojlim}$ gives the lemma. $\square$

\medskip

\noindent In other words, $\mathcal{C}(K_p,A)\leftrightarrow A[[K_p]]$ under the correspondence discussed above.

\subsection{Local Iwasawa algebras of pro-$p$-groups}

A local ring is a (possibly non-commutative) ring $R$, whose Jacobson radical $J(R)$ is a two-sided maximal ideal $\m_R$. In other words, there is a unique maximal left ideal, and a unique maximal right ideal, and they coincide. Nakayama's lemma even holds for non-commutative local rings, as is easily checked. In particular, a finitely generated projective $R$-module is free; a key fact we will make use of below, by taking $R$ to be the Iwasawa algebra of a pro-$p$-group, which turns out to be local. We first assemble the following well-known facts.

\begin{lem}
Let $K_p$ be a pro-$p$-group, and let $A$ be any $p$-ring (that is, its cardinality is a finite power of $p$, such as for $A=\mathcal{O}/\varpi^s\mathcal{O}$). Then,
\begin{itemize}
\item[(1)] Let $M$ be a left $A[[K_p]]$-module, and $H \subset K_p$ an open normal subgroup. 
Then $M/I_HM$ has a nonzero $K_p$-invariant element if $M \neq I_HM$.
\item[(2)] $K_p$ acts trivially on any simple left $A[[K_p]]$-module. 
\item[(3)] $I_{K_p}\subset J(A[[K_p]])$. 
\item[(4)] $A$ local $\Longrightarrow$ $A[[K_p]]$ local. (Furthermore, $J(A[[K_p]])=\m_A+I_{K_p}$.)
\end{itemize}
(The same is true when left modules are replaced by right modules.)
\end{lem}

\noindent {\it{Proof}}. This is all standard. We cannot resist to briefly outline the argument. For (1) it is clearly enough to show that a $p$-group $K_p$ fixes a nonzero element of any $A[K_p]$-module $M \neq 0$. This is basic group action theory; the fact that $A$ is a $p$-ring allows to count fixed points modulo $p$. For (2), if $M$ is a 
simple left $A[[K_p]]$-module, we must have $I_HM=M$ or $I_HM=0$ for all $H$. There must be some $H$ for which $I_HM\neq M$, since $M\neq 0$ is the inverse limit of all quotients $M/I_HM$. Now (1) shows that $M^{K_p}\neq 0$. By simplicity, $K_p$ acts trivially on $M$. For (3) just use that $I_{K_p}$ is generated by elements $k-1$ with $k \in K_p$. We see from (2) that $I_{K_p}$ acts trivially on any simple left $A[[K_p]]$-module, and therefore, by the very definition of the Jacobson radical, we have the inclusion as claimed.
Now (4) is immediate from (3). Indeed any maximal left ideal of $A[[K_p]]$ must be the pull-back of $\m_A$ under the augmentation map. $\square$

\medskip

\noindent We will apply this to $A[[K_p]]$, where $A=\mathcal{O}/\varpi^s\mathcal{O}$, in conjunction with Proposition 4.

\begin{prop}
Suppose $K^p$ is sufficiently small, and let $K_p \subset G(\Q_p)$ be an open pro-$p$-group. Then there exists an integer $r>0$ such that
$$
\tilde{H}^0(K^p)^{\circ}\simeq \mathcal{C}(K_p,\mathcal{O})^r
$$
as $\mathcal{O}[K_p]$-modules. Moreover, for any maximal ideal $\m \subset \T(K^p)$, the localization $\tilde{H}^0(K^p)_{\m}^{\circ}$ sits as a topologically direct summand.
\end{prop}

\noindent {\it{Proof}}. Since $A[[K_p]]$ is local, Nakayama's lemma (and Proposition 4) tells us that $H^0(K^p,A)^{\vee}$ is a {\it{free}} $A[[K_p]]$-module, of finite rank $r_s$, say. Taking the Pontryagin dual, then yields an isomorphism of smooth $A[K_p]$-modules,
$$
H^0(K^p,\mathcal{O}/\varpi^s\mathcal{O})\simeq \mathcal{C}(K_p,\mathcal{O}/\varpi^s\mathcal{O})^{r_s}.
$$
Now, we claim that $r_s$ is in fact independent of $s>0$ (and we will just write $r$ instead of $r_s$). To see this, scale both sides of the isomorphism by $\varpi$, compare the corresponding quotients, take $H$-invariants for some $H$, and compare dimensions over $\F$. This shows that $r_s=r_1$. This allows us to take the inverse limit over $s$, to obtain an isomorphism of modules over $\mathcal{O}[K_p]$,
$$
\tilde{H}^0(K^p)^{\circ}=\underset{s}{\varprojlim} H^0(K^p,\mathcal{O}/\varpi^s\mathcal{O})\simeq  \mathcal{C}(K_p,\mathcal{O})^r.
$$
In other words, an isometry $\tilde{H}^0(K^p)\simeq \mathcal{C}(K_p,L)^r$ of Banach representations of $K_p$. Finally, we may localize at any maximal ideal $\m \subset \T(K^p)$ and realize $\tilde{H}^0(K^p)_{\m}^{\circ}$ as a (topologically) direct summand of $\mathcal{C}(K_p,\mathcal{O})^r$. $\square$

\subsection{Mahler expansions and full level at $p$}

Proposition 5 already shows that the algebraic vectors are dense in $\tilde{H}^0(K^p)$ (by employing Mahler expansions, as below). In fact, this is even true for the unit ball $\tilde{H}^0(K^p)^{\circ}$. However, we can be more precise, and prove density of the smaller set of $G(\Z_p)$-locally algebraic vectors: Those $f \in \tilde{H}^0(K^p)$ such that $\langle G(\Z_p)f \rangle$ is an algebraic representation of $G(\Z_p)$.

\begin{prop}
$\tilde{H}^0(K^p)^{G(\Z_p)-alg}$ is dense in $\tilde{H}^0(K^p)$. (Same for $\tilde{H}^0(K^p)_{\m}$.)
\end{prop}

\noindent {\it{Proof}}. Pick an open normal pro-$p$-subgroup $K_p \subset G(\Z_p)$. From Proposition 5, we have an isometry $\tilde{H}^0(K^p)\simeq \mathcal{C}(K_p,L)^r$
of Banach space representations of $K_p$. We take the topological dual space $\mathcal{L}(-,L)$ on both sides, and get  
$$
\text{$\tilde{H}^0(K^p)^{\vee}\simeq L[[K_p]]^r$, $\y$ $L[[K_p]]:=L \otimes_{\mathcal{O}}\mathcal{O}[[K_p]]$.}
$$
Here $L[[K_p]]$ is identified with the distribution algebra $\mathcal{C}(K_p,L)^{\vee}$ (equipped with the bounded-weak topology) as in [ScT]. Thus, 
$\tilde{H}^0(K^p)^{\vee}$ is a free $L[[K_p]]$-module of rank $r$. It follows that $\tilde{H}^0(K^p)^{\vee}$ is projective over $L[[G(\Z_p)]]$, as
$$
\Hom_{L[[G(\Z_p)]]}(\tilde{H}^0(K^p)^{\vee},-)=\Hom_{L[[K_p]]}(\tilde{H}^0(K^p)^{\vee},-)^{G(\Z_p)/K_p}
$$
is exact: $\tilde{H}^0(K^p)^{\vee}$ is projective over $L[[K_p]]$, and taking invariants under the finite group $G(\Z_p)/K_p$ is exact, by averaging (we are in characteristic zero). 
Being projective, $\tilde{H}^0(K^p)^{\vee}$ is a direct summand of a free module (of finite rank by finite generation). That is, there is an $s>0$, and a submodule $Z$, such that
$$
\tilde{H}^0(K^p)^{\vee}\oplus Z \simeq L[[G(\Z_p)]]^s.
$$
Again, undoing the dual, and invoking Corollary 2.2 and Theorem 3.5 in [ScT],
$$
\tilde{H}^0(K^p) \oplus Z^{\vee} \simeq \mathcal{C}(G(\Z_p),L)^s.
$$
Comparing the $G(\Z_p)$-algebraic vectors on both sides, we see that it suffices to show they are dense in $\mathcal{C}(G(\Z_p),L)$. Now, topologically, we identify $G(\Z_p)\simeq \prod_{v|p}\GL_n(\mathcal{O}_{\tilde{v}})$ with a closed-open subset of $\prod_{v|p}\mathcal{O}_{\tilde{v}}^{n^2}\simeq \Z_p^t$, where we have introduced $t=[F:\Q]n^2$. Any continuous function on $G(\Z_p)$ therefore extends (non-uniquely) to a continuous function on $\Z_p^t$, which has a (multi-variable) Mahler power series expansion [Mah], which shows that the polynomials are dense in $\mathcal{C}(\Z_p^t,L)$. Finally, observe that polynomials obviously restrict to $G(\Z_p)$-algebraic functions in $\mathcal{C}(G(\Z_p),L)$. At last, localize at $\m$. $\square$

\subsection{Zariski density of crystalline points}

Following Emerton, in section 5.4 of [Eme], we deduce from the previous Proposition that "crystalline points are dense". 
\begin{cor}
The submodule $\oplus_{\lambda \in C}\tilde{H}^0(K^p)^{alg}[\lambda]$ is dense in $\tilde{H}^0(K^p)$, where $C$ denotes the collection of Hecke eigensystems $\lambda: \mathcal{H}(K^p)^{sph} \rightarrow \bar{\Q}_p$ associated with an automorphic $\pi$, which is unramified at $p$ (and of tame level $K^p$). Thus, the set of points $\ker(\lambda)$,
with $\lambda \in C$, are Zariski dense in $\text{Spec} \T(K^p)[\frac{1}{p}]$.
\end{cor}

\noindent {\it{Proof}}. First off, recall from section 3.2 that we have a decomposition,
$$
\tilde{H}^0(K^p)^{alg}=\oplus_W W \otimes H^0(K^p,\mathcal{V}_{\breve{W}})=\oplus_W \oplus_{\pi: \pi_{\infty}\simeq W}m_G(\pi)(W\otimes \pi_f^{K^p}).
$$
In particular,
$$
\tilde{H}^0(K^p)^{G(\Z_p)-alg}=\oplus_W \oplus_{\pi: \pi_{\infty}\simeq W}m_G(\pi)(W\otimes \pi_p^{G(\Z_p)}\otimes (\pi_f^p)^{K^p}),
$$
which is dense in $\tilde{H}^0(K^p)$. A fortiori, so is the $G(\Q_p)$-submodule it generates,
$$
\langle \tilde{H}^0(K^p)^{G(\Z_p)-alg}\rangle_{G(\Q_p)}=\oplus_W \oplus_{\pi: \pi_{\infty}\simeq W, \pi_p^{G(\Z_p)}\neq 0} m_G(\pi)(W\otimes \pi_f^{K^p}).
$$
We decompose the latter into eigenspaces for the action $\mathcal{H}(K^p)^{sph}$. That is, as
$$
\langle \tilde{H}^0(K^p)^{G(\Z_p)-alg}\rangle_{G(\Q_p)}=\oplus_{\lambda}\tilde{H}^0(K^p)^{alg}[\lambda]
$$
where $\lambda: \mathcal{H}(K^p)^{sph} \rightarrow \bar{\Q}_p$ runs over all eigensystems of the form $\lambda=\lambda_{\pi}$, for some automorphic $\pi$, of tamel level $K^p$, which is {\it{unramified}} at $p$ (and of some weight $W$). Thus, elements of $\cap_{\lambda\in C}\ker(\lambda)$ act trivially on $\tilde{H}^0(K^p)$.
$\square$

\section{Reduction mod $p$ and refined Serre weights}

We fix a tame level $K^p$ and point out how the eigenvariety $\X$ defines a (semisimple) global mod $p$ Langlands correspondence $t_x \leftrightarrow b_x$, analogous to the $p$-adic case discussed above. There is a natural notion of {\it{refined}} Serre weights of $\bar{\rho}_x$, and we relate them to the mod $p$ representation theory of $b_x$, and to the existence of crystalline lifts from $\X$ of compatible type $(\xi,\zeta)$.

\subsection{Integral models at $p$}

We keep the notation from previous chapters. Thus, recall that $G=\text{Res}_{F|\Q}(U)$, where $U=U(D,\star)$ is a unitary group in $n$ variables over $F$, which becomes $D^{\times}$ over $\K$. Also, we fix a prime $p$ which is {\it{unramified}} in $F$, and such that every $v|p$ of $F$ splits in $\K$. Moreover, at each $w|v$ of $\K$, we assume $D_w \simeq M_n(\K_w)$. Once and for all, above each $v|p$, we fix a place $\tilde{v}$ of $\K$, and use it to identify
$$
G(\Q_p)\simeq \prod_{v|p}\GL_n(\K_{\tilde{v}}) \hookrightarrow G(\bar{\Q}_p)\simeq \prod_{v|p}\GL_n(\bar{\Q}_p)^{\Hom(\K_{\tilde{v}},\bar{\Q}_p)}.
$$
By our assumptions, $G$ is unramified over $\Q_p$. Indeed it has a Borel pair $(B,T)$ defined over $\Q_p$, and $G$ splits over an unramified extension (the compositum of all $\K_{\tilde{v}}=F_v$). By Bruhat-Tits theory, there is a smooth affine integral model $G_{/\Z_p}$ with connected reductive special fiber. Hence $G(\Z_p)$ is hyperspecial, and
$$
G(\F_p)\simeq \prod_{v|p}\GL_n(\F_{\tilde{v}}) \hookrightarrow G(\bar{\F}_p)\simeq \prod_{v|p}\GL_n(\bar{\F}_p)^{\Hom(\F_{\tilde{v}},\bar{\F}_p)}.
$$
(Here we have tacitly used the natural bijection $\Hom(\K_{\tilde{v}},\bar{\Q}_p) \overset{\sim}{\longrightarrow} \Hom(\F_{\tilde{v}},\bar{\F}_p)$, using that $p$ is unramified.)
Moreover, one can spread out the Borel pair over $\Z_p$, and reduce it mod $p$. In particular, this allows us to identify dominant weights in different characteristics,
$$
X^*(T_{\bar{\Q}_p})_+ \simeq X^*(T_{\bar{\Z}_p})_+ \simeq X^*(T_{\bar{\F}_p})_+,
$$
with tuples $a=(a_{\tau})_{\tau \in \Hom(\K,\C)}$, where each $a_{\tau}=(a_{\tau,j})$ is a decreasing $n$-tuple of integers, such that the following polarization condition is satisfied,
$a_{\tau c,j}=-a_{\tau,n+1-j}$. This last step requires the choice of an $\iota: \C \overset{\sim}{\longrightarrow}\bar{\Q}_p$.

\subsection{Refined Serre weights}

A {\it{Serre weight}} is commonly defined as an irreducible representation $\omega$ of $G(\F_p)$ with coefficients in $\bar{\F}_p$ (usually inflated to a representation of $G(\Z_p)$ by composing with the reduction map). Since $G$ has a simply connected derived group, a result of Steinberg (extended to reductive groups by Herzig) shows that such $\omega$ are restrictions of algebraic representations $\omega$ of $G(\bar{\F}_p)$. For the restriction $\omega|_{G(\F_p)}$ to remain irreducible, the highest weight of $\omega$ has to be $p$-{\it{restricted}}. In the notation introduced in the previous subsection, this means all gaps are in the range $a_{\tau,j}-a_{\tau,j+1}<p$ (and non-negative).

\medskip

\noindent In what follows, we will simply write $G=G(\Q_p)$ and $K=G(\Z_p)$ when there is no risk of confusion. For each Serre weight $\omega$, we introduce the $\omega$-spherical Hecke algebra:
$$
\mathcal{H}_{\omega}(G,K)=\End_{G}(\text{c-Ind}_K^G(\omega)).
$$
By Frobenius reciprocity, this can be thought of more concretely as the convolution algebra of compactly supported $K$-biequivariant functions $G\rightarrow \End(\omega)$.
$\mathcal{H}_{\omega}(G,K)$ is a {\it{commutative}} noetherian $\bar{\F}_p$-algebra, according to Corollary 1.3 in [Her]. In fact, there is a mod $p$ analogue of the Satake isomorphism.

\begin{df}
A refined Serre weight is a pair $(\omega,\nu)$, consisting of a Serre weight $\omega$, together with an algebra homomorphism $
\nu: \mathcal{H}_{\omega}(G,K)\rightarrow \bar{\F}_p$.
\end{df}

\noindent Each such pair $(\omega,\nu)$ defines an $\bar{\F}_p$-representation of $G$, as the quotient
$$
\pi(\omega,\nu)=\text{c-Ind}_K^G(\omega)\otimes_{ \mathcal{H}_{\omega}(G,K),\nu}\bar{\F}_p=(\text{c-Ind}_K^G(\omega))/(\ker\nu)(\text{c-Ind}_K^G(\omega)).
$$
(Analogous to the universal modules $H_{\xi,\zeta}$ from 6.2.) If $b$ is an admissible mod $p$ representation of $G$, there is a natural action of the Hecke algebra $\mathcal{H}_{\omega}(G,K)$ on the finite-dimensional multiplicity space,
$$
\Hom_K(\omega,b)\overset{\sim}{\longrightarrow}\Hom_G(\text{c-Ind}_K^G(\omega),b),
$$
which therefore decomposes as a direct sum of generalized eigenspaces. 

\begin{lem}
The character $\nu: \mathcal{H}_{\omega}(G,K)\rightarrow \bar{\F}_p$ occurs as an eigensystem in $\Hom_K(\omega,b)$ if and only if there is a nonzero $G$-equivariant map
$\pi(\omega,\nu)\rightarrow b$.
\end{lem}

\noindent {\it{Proof}}. Immediate from the definitions. $\square$

\subsection{Types compatible with a refined Serre weight}

For the moment, we fix a refined Serre weight $(\omega,\nu)$ as above. We will define what it means for a pair $(\xi,\zeta)$ to be compatible with $(\omega,\nu)$. As in previous sections, $\xi$ denotes an irreducible algebraic $\bar{\Q}_p$-representation of $G(\bar{\Q}_p)$, and 
$$
\zeta:\mathcal{H}_{\xi}(G,K)\simeq \mathcal{H}(G,K)\longrightarrow \bar{\Q}_p
$$ 
is an algebra homomorphism. Eventually we will take $(\xi,\zeta)$ to be the {\it{type}} of a crystalline representation, in the sense used previously in this paper (that is, $\xi$ is defined from the Hodge-Tate weights, and $\zeta$ gives the eigensystem of the associated unramified representation).

\begin{df}
We say $(\xi,\zeta)$ is compatible with $(\omega,\nu)$ if there is a $G(\Z_p)$-invariant norm $\|\cdot\|_{\xi}$ (equivalently, a stable lattice $\xi^{\circ}$, the unit ball) such that
\begin{itemize}
\item $\omega \hookrightarrow \xi^{\circ}\otimes \bar{\F}_p$ (as a $G(\Z_p)$-submodule),
\item $\zeta$ is $\bar{\Z}_p$-valued on $\mathcal{H}_{\xi}(G,K)^{\circ}$, and its reduction $\zeta\otimes 1\leftrightarrow \nu$ under
$$
\mathcal{H}_{\xi}(G,K)^{\circ}\otimes \bar{\F}_p \overset{\sim}{\longrightarrow} \mathcal{H}_{\omega}(G,K).
$$
(This is the comparison isomorphism $\alpha$ from Proposition 2.10 in [Her].)
\end{itemize}
\end{df}

\noindent A few words of elaboration: The first bullet is really just saying that $\omega$ and $\xi$ have the {\it{same}} highest weight $a$. Indeed, in the notation of [Jan], if $\omega=L(a)$, we may take $\xi=H^0(a)$, a dual Weyl module, which can be defined over $\bar{\Z}_p$. However, note that for $\omega \hookrightarrow \xi^{\circ}\otimes \bar{\F}_p$ to be an isomorphism, the weight has to be $p$-{\it{small}} (according to Corollary 5.6 on p. 221 in [Jan]), which means
$$
0 \leq a_{\tau,1}-a_{\tau,n}\leq p-(n-1),
$$
for all $\tau$. (Indeed, the highest weight paired with any positive coroot $\alpha^{\vee}$ should be at most $p-\langle\varrho,\alpha^{\vee}\rangle$, where $\varrho=\frac{1}{2}\sum_{\alpha>0}\alpha$. In the $\GL_n$-case this translates into the inequalities $a_{\tau,i}-a_{\tau,j}\leq p-(j-i)$ whenever $i<j$, the strongest of which is the bound on the total gap, $a_{\tau,1}-a_{\tau,n}\leq p-(n-1)$, which is $p$-smallness.)

\medskip

\noindent Note that $p$-smallness is stronger than $p$-restrictedness for $n>2$. 

\medskip

\noindent The isomorphism in the second bullet is defined as follows: Take an $h \in \mathcal{H}_{\xi}(G,K)$, of sup-norm at most one. In particular, $h(g)$ maps 
$\xi^{\circ}$ to itself for any $g \in G$. What is shown in Proposition 2.1 in [Her] is that the reduction $\overline{h(g)}$ in fact preserves any submodule of $\xi^{\circ}\otimes \bar{\F}_p$, hence $\omega$. This yields $g \mapsto \overline{h(g)}|_{\omega}$.

\subsection{A global mod $p$ Langlands correspondence}

There is a natural mod $p$ analogue of the $p$-adic correspondence preceding Theorem B, obtained as follows. Take a point $x \in \X(L)$, where $L|E$ is an arbitrary finite extension. On the one hand, it gives rise to a pseudo-representation $t_x=\tr \bar{\rho}_x^{ss}$ (by reducing $\mathcal{T}_x=\tr \rho_x$ mod $\p_L$). On the other, we may reduce the eigensystem $\lambda_x: \frak{h}^{\circ}\rightarrow \mathcal{O}_L$ modulo $\p_L$, and look at the {\it{generalized}} eigenspace
$$
\text{$b_x=H^0(K^p,\F_L)^{ \frak{h}^{\circ}=\bar{\lambda}_x}=H^0(K^p,\F_L)_{\frak{m}_x}$, $\y$ $\frak{m}_x=\ker(\bar{\lambda}_x)$.}
$$
Here $H^0(K^p,\F_L)$ is nothing but the space of smooth functions $Y(K^p)\rightarrow \F_L$. Clearly $b_x$ carries an $\mathcal{H}_{\F_L}(K^p)$-module structure, and a commuting {\it{admissible}} action of $G(\Q_p)$. As in the $p$-adic case, the pairs $(t_x,b_x)$ form the graph of a bijection,
$$
\left\{ \begin{matrix} \text{$n$-dimensional pseudo-representations} \\ \text{$t: \Gamma_{\K,\Sigma}\rightarrow \F_L$ coming from $\X(L)$} \end{matrix} \right\}
\longleftrightarrow
$$
$$
\left\{ \begin{matrix} \text{$\mathcal{H}_{\F_L}(K^p)$-modules $b$ with admissible} \\ \text{$G(\Q_p)$-action, coming from $\X(L)$}  \end{matrix} \right\}.
$$
Here $t \leftrightarrow b$ if there is a point $x \in \X(L)$ such that $t=t_x$ and $b=b_x$, in which case
$$
t_x(\text{Frob}_w)=\bar{\lambda}_x(b_{w|v}(h_w))
$$
for all $w \notin \Sigma$. This is compatible with the $p$-adic correspondence in the following sense: For any $x' \in \X(L)$ with eigensystem $\lambda_{x'}\equiv \lambda_x$ (mod $\p_L$) we have
$$
\bar{\mathcal{B}}_{x'}=\mathcal{B}_{x'}^{\circ}\otimes \F_L \hookrightarrow b_x,
$$
a $G(\Q_p)\times \mathcal{H}_{\F_L}(K^p)$-equivariant embedding.

\subsection{Algebraic modular forms mod $p$}

We will only introduce mod $p$ modular forms of level $G(\Z_p)K^p$, and weight $\omega$,
$$
\mathcal{A}_{\omega}(G(\Z_p)K^p;\bar{\F}_p)=\Hom_{G(\Z_p)}(\omega, H^0(K^p,\bar{\F}_p)).
$$
More concretely, they are $K^p$-invariant functions $f: G(\Q)\backslash G(\A_f)\rightarrow \omega^{\vee}$ with
$$
\text{$f(gk)=\omega^{\vee}(k)^{-1}f(g)$, $\y$ $k \in G(\Z_p)$.}
$$
There is a natural action of the algebras $\mathcal{H}_{\omega}(G,K)$ and $\mathcal{H}_{\F_L}(K^p)$, and hence it decomposes as a direct sum of generalized eigenspaces for the action of the commutative algebra $\mathcal{H}_{\omega}(G,K)\otimes \mathcal{H}_{\F_L}(K^p)^{\text{sph}}$. Which eigensystems occur?

\medskip

\noindent We fix a point $x \in \X(L)$ such that $\bar{\rho}_x$ is absolutely {\it{irreducible}}, as a representation of $\Gamma_{\K}$. The Serre weights $\mathcal{W}(\bar{\rho}_x)$, relative to $K^p$, are those $\omega$ such that $\bar{\lambda}_x$ occurs in the mod $p$ modular forms. That is to say, such that
$$
\Hom_{G(\Z_p)}(\omega,b_x)\neq 0.
$$
In other words, $\text{soc}_{G(\Z_p)}(b_x)$ is of the form $\oplus_{\omega \in \mathcal{W}(\bar{\rho}_x)}m_x(\omega)\omega$. This naturally leads up to a more refined notion, taking into account the Hecke-action at $p$.

\begin{df}
$\mathcal{W}_+(\bar{\rho}_x)$ is the set of refined Serre weights $(\omega,\nu)$ such that $\nu\otimes \bar{\lambda}_x$ occurs as an eigensystem in 
the space of mod $p$ modular forms $\mathcal{A}_{\omega}(G(\Z_p)K^p;\bar{\F}_p)$.
\end{df}

\noindent Equivalently, by Lemma 1, 
$$
(\omega,\nu)\in \mathcal{W}_+(\bar{\rho}_x)\Longleftrightarrow\Hom_{G(\Q_p)}(\pi(\omega,\nu),b_x)\neq 0.
$$
(Note that $\mathcal{W}_+(\bar{\rho}_x)$ may depend on the tamel level $K^p$, at least a priori.)

\subsection{Crystalline lifts of compatible type}

The next result is inspired by Gee's approach to Serre weights via {\it{local}} crystalline lifts with prescribed Hodge-Tate weights. See [Gee], for example. We will instead consider {\it{global}} lifts, which show up on our fixed eigenvariety, and prescribe both the Hodge-Tate weights {\it{and}} the Frobenius eigenvalues. 

\medskip

\noindent For technical reasons, we will assume from now on that our fixed tame level $K^p$ is sufficiently small (for precision, that some $K_v$ has no elements of order $p$).

\begin{prop}
Let $(\omega,\nu)$ be an arbitrary refined Serre weight, and let $x$ be a point of $X_{reg}$ such that $\bar{\rho}_x$ is {\it{irreducible}}. Then the following 
holds.
\begin{itemize}
\item[(a)] If $(\omega,\nu)\in \mathcal{W}_+(\bar{\rho}_x)$, then $\bar{\rho}_x$ has a crystalline lift from $X_{reg}$, whose type is compatible with $(\omega,\nu)$.
\item[(b)] The converse holds when $\omega$ is $p$-\underline{small}: If $\bar{\rho}_x$ has a crystalline lift from $X_{reg}$, whose type is compatible with $(\omega,\nu)$, then $(\omega,\nu)\in \mathcal{W}_+(\bar{\rho}_x)$.
\end{itemize}
\end{prop}

\noindent {\it{Proof of (a)}}:

\medskip

\noindent Say $\omega=L(a)$. Pick $\xi^{\circ}$ to be a $\bar{\Z}_p$-structure in the dual Weyl module $\xi=H^0(a)$. Then $\omega$ is the unique irreducible submodule of $\xi^{\circ}\otimes \bar{\F}_p$. This yields a natural reduction map,
$$
\mathcal{A}_{\xi}(G(\Z_p)K^p;\bar{\Z}_p)\otimes \bar{\F}_p \rightarrow \mathcal{A}_{\omega}(G(\Z_p)K^p;\bar{\F}_p),
$$
where
$$
\mathcal{A}_{\xi}(G(\Z_p)K^p;\bar{\Z}_p)=\Hom_{G(\Z_p)}(\xi^{\circ}, H^0(K^p,\bar{\Z}_p))
$$
is naturally a $\bar{\Z}_p$-structure in a space of classical $p$-adic modular forms. The reduction map is easily seen to be surjective if $K^p$ is sufficiently small (that is, if $p$ does not divide the orders of the arithmetic subgroups showing up). By assumption, the inflated eigensystem $\nu\otimes \bar{\lambda}_x$,
$$
\mathcal{H}_{\xi}(G,K)^{\circ}\otimes \frak{h}^{\circ} \longrightarrow \mathcal{H}_{\omega}(G,K)\otimes \mathcal{H}_{\F_L}(K^p)^{\text{sph}}\longrightarrow \bar{\F}_p,
$$
occurs in $\mathcal{A}_{\omega}(G(\Z_p)K^p;\bar{\F}_p)$. Localizing the reduction map at $\ker(\nu\otimes \bar{\lambda}_x)$, and invoking the Deligne-Serre lifting lemma, we infer that there is an eigensystem 
$$
\Lambda: \mathcal{H}_{\xi}(G,K)^{\circ}\otimes \frak{h}^{\circ} \longrightarrow \bar{\Z}_p,
$$
with $\bar{\Lambda}=\nu\otimes \bar{\lambda}_x$, which occurs in the space $\mathcal{A}_{\xi}(G(\Z_p)K^p;\bar{\Q}_p)$ of classical $p$-adic modular forms. Via 
$\iota: \C \overset{\sim}{\longrightarrow}\bar{\Q}_p$ we thus find an automorphic representation $\pi$ of $G(\A)$ with $\pi_{\infty}=\xi$, such that $\frak{h}$ acts on $\pi_f^{K^p}\neq 0$ by the character $\Lambda^p\equiv \lambda_x$, and such that $\pi_p$ is unramified with eigensystem $\Lambda_p$ lifting $\nu$ (where we use the usual identification 
$\mathcal{H}(G,K)\simeq \mathcal{H}_{\xi}(G,K)$ to view $\Lambda_p$ as a character of the spherical Hecke algebra). Say $\pi_p\subset \text{Ind}(\theta)$, and $\psi$ is the highest weight of $\xi$. Then $(\psi\theta,\Lambda^p)$ belongs to $E(0,K^p)_{cl}$. Let $x' \in X_{cl}$ be the corresponding classical point of the eigenvariety. Now, the Galois representation $\rho_{x'}\simeq \rho_{\pi,\iota}$ (comes from $X_{cl}$ and) is crystalline at all places avove $p$; since $\pi_p$ is unramified. Moreover, $\bar{\rho}_{x'}\simeq \bar{\rho}_x$ since $\Lambda^p\equiv \lambda_x$ (both are irreducible). Finally, $\rho_{x'}$ has type $(\xi,\Lambda_p)$, which is compatible with $(\omega,\nu)$. Indeed
$\omega \subset \xi^{\circ}\otimes \bar{\F}_p$, and obviously $\bar{\Lambda}_p\leftrightarrow \nu$.

\medskip

\noindent {\it{Proof of (b)}}:

\medskip

\noindent Now suppose there is an old point $x' \in X_{reg}$ such that $\bar{\rho}_{x'}\simeq \bar{\rho}_x$, and $\rho_{x'}$ has type $(\xi,\zeta)$ compatible with $(\omega,\nu)$.
By Theorem 4, part (3), we know that
$$
\xi_{x'}\otimes \text{Ind}_B^G(\theta_{x'})=\widetilde{BS}(\rho_{x'})\hookrightarrow \mathcal{B}_{x'},
$$
as a $G(\Q_p)$-representation. (Recall that $\rho_{x'}$ is even residually irreducible, so $\pi_{x'}$ is generic at all split places). Here $\xi=\xi_{x'}$, and $\mathcal{H}(G,K)$ acts on 
the newvectors in $\text{Ind}_B^G(\theta_{x'})$ via $\zeta$. Choosing a newvector (unique up to a constant) then gives a $K$-map
$\xi \hookrightarrow \mathcal{B}_{x'}$ (or, equivalently, a nonzero $G$-map $\text{c-Ind}_K^G(\xi)\rightarrow \mathcal{B}_{x'}$) on which $\mathcal{H}_{\xi}(G,K)\simeq \mathcal{H}(G,K)$ acts by $\zeta$. Let $\|\cdot\|_{\xi}$ be the norm on $\xi$ (and $\xi^{\circ}$ its unit ball), whose existence is guaranteed by Definition 9; the compatibility with $(\omega,\nu)$.
Since $\xi$ is finite-dimensional, all norms are equivalent. So $\xi \hookrightarrow \mathcal{B}_{x'}$ is automatically continuous, and upon scaling, we may assume it preserves the unit balls. That is, $\xi^{\circ}\hookrightarrow \xi \cap \mathcal{B}_{x'}^{\circ}$, with cokernel killed by some nonzero constant in $\bar{\Z}_p$ (all lattices are commensurable). 
By the Brauer-Nesbitt principle, 
$$
(\xi^{\circ}\otimes \bar{\F}_p)^{ss}\simeq ((\xi \cap \mathcal{B}_{x'}^{\circ}) \otimes \bar{\F}_p)^{ss}.
$$
Moreover, the cokernel of the inclusion $\xi \cap \mathcal{B}_{x'}^{\circ}\hookrightarrow \mathcal{B}_{x'}^{\circ}$ is clearly torsion-free, so
$$
(\xi \cap \mathcal{B}_{x'}^{\circ}) \otimes \bar{\F}_p\hookrightarrow \mathcal{B}_{x'}^{\circ} \otimes \bar{\F}_p\hookrightarrow b_x,
$$
since $\lambda_{x'}\equiv \lambda_{x}$. Now, since we assume $\omega$ is $p$-small, the source of this map is isomorphic to $\omega$, and the way it was constructed shows that $\mathcal{H}_{\omega}(G,K)$ acts on $\omega \hookrightarrow b_x$ by $\nu\leftrightarrow \zeta\otimes 1$. Ths results in a nonzero $G$-map $\pi(\omega,\nu)\rightarrow b_x$. $\square$

\medskip

\noindent {\it{Remark}}. If we are not assuming $p$-smallness in (b), the above argument still shows that $\mathcal{W}(\bar{\rho}_x)$ contains {\it{some}} constituent of 
$\xi^{\circ}\otimes \bar{\F}_p$ (of highest weight at {\it{most}} that of $\omega$), but we cannot say much about the Hecke action at $p$.

\section{Local-global compatibility and Ihara's lemma}

We define a candidate $b(\bar{\rho}_x)$ for $\otimes_{v|p}b(\bar{\rho}_x|_{\Gamma_{\K_{\tilde{v}}}})$ by means of the local Langlands correspondence in characteristic $p$ (see Theorem 5.1.5 in [EH]). We offer some (precise) {\it{speculations}} towards local-global compatibility mod $p$, from which the conjectural Ihara lemma (Conjecture B in [CHT]) is deduced.

\subsection{Local Langlands in characteristic $p$}

For finite extensions $K|\Q_{\ell}$, where $\ell \neq p$, Emerton and Helm have defined a map
$$
\left\{ \begin{matrix} \text{continuous $\bar{\rho}:\Gamma_K \rightarrow \GL_n(\F_L)$}  \end{matrix} \right\}
\longrightarrow
$$
$$
\left\{ \begin{matrix} \text{$\F_L$-spaces $\bar{\pi}$ endowed with} \\ \text{admissible $\GL_n(K)$-action}  \end{matrix} \right\},
$$
denoted $\bar{\rho}\mapsto \bar{\pi}(\bar{\rho})$. This $\bar{\pi}(\bar{\rho})$ is of finite length, and is characterized uniquely by the following three properties (according to
Theorem 5.1.5 in [EH]): 
\begin{itemize}
\item[(a)] $\bar{\pi}(\bar{\rho})$ is essentially AIG.
\item[(b)] If $\rho: \Gamma_K \rightarrow \GL_n(\mathcal{O}_{L'})$ is a continuous lift of $\bar{\rho}\otimes \F_{L'}$, and $\Lambda$ is a $\GL_n(K)$-invariant lattice in $\pi(\rho)$ (the generic local Langlands correspondent) with $\Lambda \otimes \F_{L'}$ essentially AIG, then there is a $\GL_n(K)$-embedding,
$$
\Lambda \otimes \F_{L'} \hookrightarrow \bar{\pi}(\bar{\rho})\otimes \F_{L'}.
$$
\item[(c)] $\bar{\pi}(\bar{\rho})$ is minimal in the following sense: If $\bar{\pi}$ is any representation satisfying (a) and (b), then there is a $\GL_n(K)$-embedding, $\bar{\pi}(\bar{\rho}) \hookrightarrow \bar{\pi}$.
\end{itemize}
The correspondence has a number of additional nice properties, which we will not recall. See (4) through (8) on p. 44 in [EH].

\medskip

\noindent AIG stands for absolutely irreducible and generic. Paraphrasing Definition 3.2.1 in [EH], a smooth representation $\bar{\pi}$ is {\it{essentially AIG}} if
$\text{soc}(\bar{\pi})$ is AIG, $\bar{\pi}/\text{soc}(\bar{\pi})$ has no generic constituents, and $\bar{\pi}$ is the sum of its finite length submodules.

\medskip

\noindent The correspondence $\bar{\rho}\mapsto \bar{\pi}(\bar{\rho})$ plays a key role in (strong) local-global compatibility mod $p$ for $\GL(2)$. See Theorem 1.2.6 in [Eme], which will serve as a guide in what follows.

\subsection{Conjectural local-global compatibility mod $p$}

We now return to the situation from Chapter 6. Thus, we have a tame level $K^p$ with {\it{split ramification}}, and we let $x \in \X(L)$ be a regular (classical) point on the eigenvariety such that $\rho_x$ is {\it{residually}} irreducible, and $m(\pi_x)=1$.

\begin{df}
For each place $v \nmid p$ of $F$, we define a representation $\bar{\pi}_{x,v}$ of $U(F_v)$ over $\F_L$. In the unramified case, where $v \notin S(K^p)$, $\pi_{x,v}$ is a $p$-integral unramified principal series, and we define $\bar{\pi}_{x,v}$ to be its canonical reduction. In the ramified case, where $v \in S(K^p)$ necessarily splits in $\K$, 
$$
\bar{\pi}_{x,v}\simeq \bar{\pi}(\bar{\rho}_x|_{\Gamma_{\K_w}}),
$$
under $U(F_v)\simeq \GL_n(\K_w)$, where $w|v$ is any of the two places above.
\end{df}

\noindent Only the $\mathcal{H}(K^p)$-module $m=\otimes_{v\nmid p}\bar{\pi}_{x,v}^{K_v}$ will occur below, on which $\mathcal{H}(K^p)^{\text{sph}}$ acts by $\bar{\lambda}_x$.
Thus the precise definition of $\bar{\pi}_{x,v}$ at $v \notin S(K^p)$ is not relevant here. The ramified places $S(K^p)$ are the interesting ones, where we invoke [EH].

\medskip

\noindent Guided by Theorem 1.2.6 in [Eme], and the remark after Corollary 3 above, we define our candidate representation $b(\bar{\rho}_x)$ as follows.

\begin{df}
$b(\bar{\rho}_x)=\Hom_{\mathcal{H}(K^p)}(\otimes_{v\nmid p}\bar{\pi}_{x,v}^{K_v},b_x)$.
\end{df}

\noindent This is an admissible representation of $G(\Q_p)$ over $\F_L$. However, it is not clear at all whether $b(\bar{\rho}_x)$ is {\it{nonzero}}! This is implied by local-global compatibility:

\begin{cn}
The tautological map
$$
\Phi:b(\bar{\rho}_x) \otimes (\bigotimes_{v\nmid p}\bar{\pi}_{x,v}^{K_v}) \rightarrow b_x
$$ 
is an isomorphism (preserving the action of $G(\Q_p)$ and $\mathcal{H}(K^p)$).
\end{cn}

\medskip

\noindent This may very well turn out to be too naive. For now, it serves as a guiding principle. In the $p$-adic case, we had to "replace" $\mathcal{B}_x$ with the closure of the regular-algebraic vectors $\overline{\mathcal{B}_x^{ralg}}$ to achieve local-global compatibility on the nose. Maybe one should look for an analogous (finite length?)
submodule $b_x^{\#}\subset b_x$, containing all the reductions $\bar{\mathcal{B}}_{x'}\hookrightarrow b_x^{\#}$, and use {\it{it}} to define $b(\bar{\rho}_x)$ as above.

\medskip

\noindent Of course, $b(\bar{\rho}_x)$ appears to be nothing more than a toy construction unless one can show it only depends on the various restrictions
$\bar{\rho}_x|_{\Gamma_{\K_{\tilde{v}}}}$ at $p$, and that it factors as a tensor product $\otimes_{v|p}$. This is clearly at the heart of the whole discussion! As in the $p$-adic case, we have made no progress towards this problem.

\medskip

\noindent {\it{Remark}}. Conjecture 1 holds if $S(K^p)=\varnothing$, in which case $b(\bar{\rho}_x)\simeq b_x$.

\subsection{Conjectural $p$-adic and mod $p$ compatibility}

We keep the assumptions and the notation of the previous section. Thus, we fix  a point $x \in \X(L)$ such that $\bar{\rho}_x$ is absolutely irreducible (as a representation of $\Gamma_{\K}$). There should be a relation between $b(\bar{\rho}_x)$ and the reduction of $B(\rho_x)$, the $p$-adic candidate introduced in Theorem B of the introduction (and in the remark after Corollary 3 in section 6.1). We will make this expectation precise.

\medskip

\noindent First, we will specify a lattice in $M=\otimes_{v\nmid p}\pi_{x,v}^{K_v}$. Let $v\nmid p$ be a ramified (hence split) place of $F$. By Proposition 3.3.2 in [EH], since $\pi_{x,v}$ is generic, it admits a $U(F_v)$-invariant lattice $\pi_{x,v}^{\circ}$ such that its (naive) reduction $\pi_{x,v}^{\circ}\otimes \F_L$ is essentially AIG (which is unique up to homothety). By (2) of Theorem 5.1.5 in [EH], one can then find a $U(F_v)$-equivariant embedding $\pi_{x,v}^{\circ}\otimes \F_L\hookrightarrow \bar{\pi}_{x,v}$ (not necessarily onto, when the target is reducible). As lattice in $M$, we take $M^{\circ}=\otimes_{v\nmid p} (\pi_{x,v}^{\circ})^{K_v}$.

\begin{cn}
There is a $G(\Q_p)$-equivariant map $\Psi$ such that the diagram
$$
\xymatrix{
B(\rho_x)^{\circ}\otimes \F_L  \ar[d] \ar[r]^{\Psi} &b(\bar{\rho}_x) \ar[d]\\
\Hom_{\mathcal{H}(K^p)^{\circ}}(M^{\circ}\otimes \F_L,\mathcal{B}_x^{\circ}\otimes \F_L) \ar[r]  & \Hom_{\mathcal{H}(K^p)^{\circ}}(M^{\circ}\otimes \F_L,b_x)}
$$
commutes. (Here the right vertical map is composition with $M^{\circ}\otimes \F_L \hookrightarrow m$.)
\end{cn}

\noindent Note that such a map $\Psi$ is necessarily nonzero ($\Rightarrow b(\bar{\rho}_x)$ is nonzero).

\subsection{Connections to Ihara's lemma}

A key step in Wiles's proof of Fermat's Last Theorem was to reduce the non-minimal case of modularity lifting to the minimal case, by level-raising. What makes this work is a lemma of Ihara. This approach was mimicked in [CHT] for $U(n)$, where they deduce modularity lifting $R=\T$ from a conjectural analogue of Ihara's lemma (Conjecture B in [CHT]), originating from Mann's Harvard Ph.D. thesis, which is still open. Taylor has since adapted techniques of Kisin to bypass Ihara's lemma and prove modularity lifting theorems of the form $R^{\text{red}}=\T$. As mentioned in the introduction of [CHT], it would be interesting to get rid of the nilradical. For example, this would have applications to special values of $L$-functions.

\medskip

\noindent Let us recast Conjecture B of [CHT] in our setup:

\medskip

\noindent {\bf{Ihara's Lemma}}. {\it{Let $u \in S(K^p)$ be a place of $F$ such that $U(F_u)\simeq \GL_n(F_u)$, and $p$ is banal for this group. Let $x\in X_{cl}$ be a point such that $\bar{\rho}_x$ is irreducible.
$$
\mathcal{S}_x\overset{df}{=}\text{$\mathcal{C}^{\infty}(U(F)\backslash U(\A_{F,f})/K^{p,u},\bar{\F}_p)_{\m_x}$, $\y$ $\m_x\overset{df}{=}\ker(\bar{\lambda}_x)\subset \mathcal{H}(K^p)^{\text{sph}}$.}
$$
Then every simple $\GL_n(F_u)$-submodule of $\mathcal{S}_x$
is generic.}}

\medskip

\noindent This has close ties to Conjecture 1 (local-global compatibility mod $p$). 

\begin{prop}
Conjecture 1 $\Longrightarrow$ Ihara's lemma.
\end{prop}

\noindent {\it{Proof}}. Let $\tau \subset \mathcal{S}_x$ be a simple $\GL_n(F_u)$-submodule. In $U(F_u)$ we pick a small enough compact open subgroup $K_u$ such that $\tau^{K_u}\neq 0$, and such that taking $K_u$-invariants defines an equivalence of categories between the category of smooth representations of $U(F_u)$, which are generated by their
$K_u$-invariants, and the category of modules for the Hecke algebra $\mathcal{H}(U(F_u),K_u)$. This is possible since $p$ is assumed banal for $U(F_u)$, by a result of Vigneras. 

\medskip

\noindent Consequently, with a possibly smaller $K^p$ (namely $K_uK^{p,u}$) we have 
$$
\tau^{K_u}\subset \mathcal{S}_x^{K_u}=H^0(K^p,\bar{\F}_p)_{\m_x}=b_x.
$$
Admitting Conjecture 1, we conclude that the simple Hecke module $\tau^{K_u}$ embeds into $\bar{\pi}_{x,u}^{K_u}$. By choice of $K_u$, this arises from an embedding $\tau \hookrightarrow \bar{\pi}_{x,u}$. Hence,
$$
\tau=\text{soc}(\bar{\pi}_{x,u}),
$$
since the latter is AIG, by the desiderata in [EH]. In particular, $\tau$ is generic. $\square$

\medskip

\noindent  Admittedly, this just replaces Ihara's lemma by a stronger conjecture. However, we feel that the latter is more conceptual, and more intuitive.

\medskip

\noindent {\it{Remark}}. Conversely, it looks like something along the lines of Ihara's lemma is needed to even show that $b(\bar{\rho}_x)\neq 0$. Indeed, this would require producing maps $\bar{\pi}_{x,u}^{K_u}\rightarrow b_x$. This would make use of the minimality (c). However, to do that, one would have to verify that $\mathcal{S}_x$ is AIG, among other things, such as dealing with the lifts in (b) above. It would be interesting to see if the modularity lifting theorems in [CHT], which {\it{result}} from Ihara's lemma, shed any light on this.

\appendix
\section{Appendix: Frobenius reciprocity}

For the convenience of the reader, we briefly recall the explicit formulas giving Frobenius reciprocity, and as a consequence deduce that {\it{continuity}} is preserved. This was used in section 6 on weak local-global compatibility. 

\medskip

\noindent We will work in the following setup: Let $G$ be a topological group, and let $K$ be a compact open subgroup. Let $\xi:K \rightarrow \GL(V)$ be a continuous representation on a finite-dimensional $L$-vector space $V$, where $L/\Q_p$ is finite. The compact induction $\text{c-ind}_K^G\xi$ consists of all compactly supported functions $f:G \rightarrow V$ such that $f(kg)=\xi(k)f(g)$. Such $f$ are automatically continuous. 

\medskip

\noindent {\it{Example}}. For any vector $v \in V$ define $\phi_v$ by letting $\phi_v(g)=\xi(g)v$ when $g \in K$, and $\phi_v(g)=0$ otherwise. Observe that $\phi_{\xi(k)v}=k\phi_v$.
Consequently, these $\{\phi_v\}$ generate $\text{c-ind}_K^G\xi$ as a $G$-representation: Indeed, any $f$ as above can be written as a (finite) sum $f=\sum_{g \in K\backslash G}g^{-1}\phi_{f(g)}$.

\begin{lem}
Let $\rho:G \rightarrow \GL(W)$ be any (possibly infinite-dimensional) representation of $G$ on an $L$-vector space $W$. Then there is a natural bijection
$$
\Hom_K(\xi,\rho|_K)\overset{\sim}{\longrightarrow} \Hom_G(\text{c-ind}_K^G\xi,\rho).
$$
Choose a $K$-invariant norm $\|\cdot\|_{\xi}$ on $V$, and let $\|\cdot\|_{\xi,\infty}$ be the corresponding $G$-invariant supremum-norm on $\text{c-ind}_K^G\xi$. Moreover, suppose $W$ has a $G$-invariant norm $\|\cdot\|_{\rho}$. {\bf{Then}} the bijection preserves the transformation norms on both sides. (In particular, if $\xi$ is irreducible, every transformation $\text{c-ind}_K^G\xi\rightarrow \rho$ is continuous.)
\end{lem}

\noindent {\it{Proof}}. The two adjunction maps 
$$
\text{$\alpha: \xi \rightarrow \text{c-ind}_K^G(\xi)|_K$, $\y$ $\beta: \text{c-ind}_K^G(\rho|_K)\rightarrow \rho$,}
$$
are defined by: $\alpha(v)=\phi_v$ and $\beta(f)=\sum_{g \in K\backslash G} \rho(g)^{-1}f(g)$. Easily checked to be equivariant under $K$ and $G$ respectively. They define the bijection as follows. For $\ell \in \Hom_K(\xi,\rho|_K)$, we first induce it to a map $i(\ell): \text{c-ind}_K^G(\xi) \rightarrow \text{c-ind}_K^G(\rho|_K)$, and then compose it with $\beta$. More explicitly, with $\ell$ we associate 
$$
r(f)={\sum}_{g \in K \backslash G}\rho(g)^{-1}\ell(f(g)).
$$
Conversely, starting with an $r \in \Hom_G(\text{c-ind}_K^G\xi,\rho)$, we define an $\ell$ by first restricting $r$ to $K$, and then composing with $\alpha$. In other words, by the formula
$$
\ell(v)=r(\phi_v).
$$
It is straightforward to check these maps are mutually inverse. Suppose $\ell \leftrightarrow r$. Then, first of all, if $r$ is bounded, with transformation norm $\|r\|$, it follows that
$$
\|\ell(v)\|_{\rho}\leq \|r\|\cdot \|\phi_v\|_{\xi,\infty}=\|r\|\cdot\|v\|_{\xi},
$$
so that $\ell$ is bounded, with transformation norm $\|\ell\|\leq \|r\|$. Conversely,
$$
\|r(f)\|_{\rho}\leq {\max}_{g\in K\backslash G}\|\ell(f(g))\|_{\rho}\leq \|\ell\|\cdot \|f\|_{\xi,\infty},
$$
and hence $\|r\|\leq \|\ell\|$. Altogether, this must be an equality. If $\xi$ is irreducible, every $K$-equivariant map $\ell: \xi\rightarrow \rho|_K$ must be continuous: It is necessarily injective, if nonzero, and therefore $v \mapsto \|\ell(v)\|_{\rho}$ defines a norm on $\xi$, which must be equivalent to $\|\cdot\|_{\xi}$ since $\xi$ is finite-dimensional, see 4.13 in [Sc]. $\square$

\medskip

\noindent {\it{Remark}}. The above argument shows that $r$ isometric $\Rightarrow$ $\ell$ isometric; but {\it{not}} conversely. If $\ell$ is isometric, all we can say a priori is that
$\|r(f)\|_{\rho}=\|f\|_{\xi,\infty}$ when $f$ is of the form $f=\phi_v$, for some $v$. Or, more generally, for functions $f$ such that $g \mapsto \|f(g)\|_{\xi}$ has "no repeated
maximum", by which we mean it attains its maximum at a unique coset in $K \backslash G$.



\noindent {\sc{Department of Mathematics, Princeton University, USA.}}

\noindent {\it{E-mail address}}: {\texttt{claus@princeton.edu}}

\end{document}